\documentclass[matprg,nosmartrunhead]{svjour}

\usepackage{graphics}
\usepackage{epsf}
\usepackage{amssymb}
\def\Arg{\mbox{\rm Arg}}
\def\Var{\mbox{\rm Var}}
\def\SA{{\mathcal A}}
\def\SB{{\mathcal B}}
\def\SH{{\mathcal H}}
\def\SI{{\mathcal I}}
\def\SL{{\mathcal L}}
\def\So{{o}}
\def\SO{{\mathcal O}}
\def\SR{{\mathcal R}}
\def\SV{{\mathcal V}}
\def\SX{{\mathcal X}}
\def\Ib{{\bf I}}
\def\Mb{{\bf M}}
\def\Nb{{\bf N}}
\def\Pb{{\bf P}}
\def\Qb{{\bf Q}}

\def\me{\epsilon}
\def\mg{\gamma}
\def\TT{^\top}
\def\mt{\theta}
\def\ml{\lambda}
\def\ra{\rightarrow}
\newtheorem{coro}{\bf Corollary}
\def\essinf{\mbox{\rm ess inf}}
\def\esssup{\mbox{\rm ess sup}}
\newcommand{\carre} {\hfill \rule{2mm}{2mm} }

\begin{document}
\author{Luc Pronzato \and Henry P.\ Wynn \and Anatoly A.\ Zhigljavsky}
\title{Asymptotic behaviour of a family of gradient algorithms in $\mathbb{R}^d$
and Hilbert spaces}
\combirunning{L.\ Pronzato et. al. Asymptotic behaviour of a
family of gradient algorithms}

\institute{Luc Pronzato\at Laboratoire I3S, CNRS - UNSA, Les
Algorithmes - B\^at.\ Euclide B, 2000 route des Lucioles -- B.P.
121, F-06903 Sophia Antipolis Cedex, France. Fax: 33 (0)4 92 94 28
96; Email: pronzato@i3s.unice.fr \and Henry P.\ Wynn\at London
School of Economics, Department of Statistics, London, WC2A 2AE,
UK \and Anatoly A.\ Zhigljavsky\at Cardiff University, School of
Mathematics, Senghennydd Road, Cardiff, CF24 4AG, UK}
\date{Received: date / Revised version: date}

\maketitle
\begin{abstract}
The asymptotic behaviour of a family of gradient algorithms
(including the methods of steepest descent and minimum residues)
for the optimisation of bounded quadratic operators in
$\mathbb{R}^d$ and Hilbert spaces is analyzed. The results
obtained generalize those of Akaike (1959) in several directions.
First, all algorithms in the family are shown to have the same
asymptotic behaviour (convergence to a two-point attractor), which
implies in particular that they have similar asymptotic
convergence rates. Second, the analysis also covers the Hilbert
space case. A detailed analysis of the stability property of the
attractor is provided.
\end{abstract}
%

\section{Introduction}

The paper generalizes the results presented in \cite{PWZa01-AAM}
to other optimisation algorithms of the gradient type. We
introduce a class of algorithms, called $P$-gradient algorithms,
that differ by the choice of the length of the step made in the
gradient direction. The class includes in particular the usual
steepest-descent algorithm and the method of minimal residues of
Krasnosel'skii and Krein \cite{KozjakinK82,KrasnoselskiiK52}. We
show that for a quadratic function, the worst asymptotic rate of
convergence is the same for the whole class of algorithms
considered. It is also true that, expressed in the right
framework, all the algorithms in the class behave in a very
similar fashion\footnote{Not all algorithms using the gradient
direction belong to that class, which in particular does {\em not}
include the spectral-gradient algorithm, see \cite{BirginM2001},
proposed by Barzilai and Borwein in \cite{BarzilaiB88}. This
method, which has been found in particular examples to allow
significant improvement over standard steepest descent, see
\cite{Raydan97}, thus requires a separate treatment. The same is
true for steepest descent with relaxation or the combination of
steepest descent and Barzilai-Borwein methods, as considered in
\cite{RaydanS2002}.}. This analysis complements that presented in
\cite{Akaike59}, \cite{NocedalSZ98,NocedalSZ2002} and Chapter 7 of
\cite{PWZl2000} which concerns steepest descent. Moreover, the
analysis in \cite{PWZa01-AAM} directly applies to all algorithms
in the class considered, revealing the asymptotic behaviour for
bounded quadratic operators not only in $\mathbb{R}^d$ but also in
Hilbert spaces. The worst case behaviour exhibited is fundamental
``bottom-line'' in the study of optimisation whose understanding
is critical for building more complex and faster algorithms.


The basic idea is renormalisation, as used throughout
\cite{PWZl2000}. The main result in the finite dimension case is
that for any algorithm in the class, in the renormalised space one
observes convergence to a two-point attractor which lies in the
space spanned by the eigenvectors corresponding to the smallest
and largest eigenvalues of the matrix $A$ of the quadratic
operator. The proof for bounded quadratic operators in Hilbert
space stems from the proof for $\mathbb{R}^d$ but is considerably
more technical. In both cases, as in \cite{Akaike59}, the method
consists of converting the problem to one containing a special
type of operator on measures on the spectrum of the operator. The
additional technicalities arise from the fact that in the Hilbert
space case the measure, which is associated with the spectral
measure of the operator, may be continuous. Another important
result concerns bounds on convergence rates, named after
Kantorovich, see \cite{KantorovichA82}. For all algorithms in the
family considered, the actual asymptotic rate of convergence,
although satisfying Kantorovich bounds, depends on the starting
point and is difficult to predict. This complex behaviour has
consequences for the stability of the attractor, which are
discussed following the main results.

The family of gradient algorithms we consider, called $P$-gradient
algorithms, is introduced in Section~\ref{S:family}.
Renormalisation is presented there, which, together with the
monotonic sequences of Section~\ref{S:monotonicity}, forms the
core of the analysis to be conducted. The main results are
presented in Section~\ref{S:asymptoticbehaviour}, first for the
case $\SH=\mathbb{R}^d$, then for the Hilbert space case. They
rely on the convergence property of successive transformations of
a probability measure, which is presented in
Section~\ref{S:theorem-tansform-pm}. Again, the two cases
$\SH=\mathbb{R}^d$ and $\SH$ a Hilbert space are distinguished,
the exposition being much simpler in the former case. The
stability of attractors is discussed in Section~\ref{S:stability},
only in the more general case of a Hilbert space, the case
$\SH=\mathbb{R}^d$ not allowing for a significant simplification
of the presentation. Finally, Section~\ref{S:rate} shows the
asymptotic equivalence between several rates of convergence of
gradient algorithms. All proofs and some important lemmas are
collected in an appendix.

\section{A family of gradient algorithms}
\label{S:family}

\subsection{$P$-gradient algorithms}
\label{S:P-gradient} Let $A$ be a real bounded self-adjoint
(symmetric) operator in a real Hilbert space $\SH$ with inner
product $(x,y)$ and norm given by $\|x\|=(x,x)^{1/2}$. Assume that
$A$ is positive, bounded below, and denote its spectral boundaries
by $m$ and $M$:
$$
m = \inf_{\|x\|=1} (Ax,x)\,, \ M = \sup_{\|x\|=1} (Ax,x)\,,
$$
with $0<m<M<\infty$. The function to be minimized corresponds to
the quadratic form
\begin{equation}
f(x)= \frac{1}{2} (Ax,x) - (x,y) \,. \label{fx}
\end{equation}
It is minimum at $x^*=A^{-1}y$, its directional derivative at $x$
in the direction $u$ is
$$
\nabla_u f(x) = (Ax-y,u) \,.
$$
The direction of steepest descent at $x$ is $-g$, with $g=g(x)$
the gradient at $x$, namely $g=Ax-y$. The minimum of $f$ in this
direction is obtained for the optimum step-length
$$
\gamma=\frac{(g,g)}{(Ag,g)} \,,
$$
which corresponds to the usual steepest-descent algorithm. One
iteration of the steepest descent algorithm is thus
\begin{equation}
x_{k+1} = x_k - \frac{(g_k,g_k)}{(Ag_k,g_k)} g_k \,, \label{algo}
\end{equation}
with $g_k=Ax_k-y$ and $x_0$ some initial element in $\SH$. We
define more generally the following class of algorithms.

\begin{definition}
\label{D:P-gradient} Let $P(\cdot)$ be a real function defined on
$[m,M]$, infinitely differentiable, with Laurent series
$$
P(z) = \sum_{-\infty}^\infty c_k z^k\,, \ c_k\in\mathbb{R} \mbox{
for all } k\,,
$$
such that $0<\sum_{-\infty}^\infty c_ka^k<\infty$ for $a\in[m,M]$.
The $k$-th iteration of a $P$-gradient algorithm is defined by
\begin{equation}\label{iteration}
    x_{k+1} = x_k -\gamma_k g_k
\end{equation}
where the step-length $\gamma_k$ minimizes $(P(A)g_{k+1},g_{k+1})$
with respect to $\gamma$, with $g_{k+1}=g(x_{k+1})=g(x_k - \gamma
g_k)$.
\end{definition}

Direct calculation gives
\begin{equation}\label{gammak}
    \gamma_k= \frac{(P(A)Ag_{k},g_{k})}{(P(A)A^2g_{k},g_{k})} \,.
\end{equation}
Note that $AP(A)=P(A)A$ and that the denominator and numerator of
$\gamma_k$ are linear in $P(A)$. Also, $\gamma_k$ is
scale-invariant in $P(A)$ and $\gamma_k\in[1/M,1/m]$.

Taking $P(A)=A^{-1}$ gives the steepest-descent algorithm.
Choosing $P(A)=I$, the identity operator, is equivalent to
choosing the step-length that minimizes the norm of the gradient
$g_{k+1}$ at the next point. We then obtain the method of minimal
residues introduced in \cite{KrasnoselskiiK52} for the solution of
linear equations. For any fixed $\alpha\in(0,1)$, choosing
$\gamma_k$ that minimizes $ \alpha f(x_k - \gamma g_k)+$
$\mbox{$(1-\alpha)(g(x_k - \gamma g_k),g(x_k - \gamma g_k))$}$
with respect to $\mg$ also gives an algorithm in the family. More
generally, we show below how to construct $P$-gradient algorithms,
with $P(\cdot)$ a polynomial in $A$, using evaluations of
$f(\cdot)$ and $g(\cdot)$ only.

\subsection{Practical construction when $P$ is a polynomial}
\label{S:constructions} We consider the case where $P(A)=A^q$ for
some integer $q\geq -1$. (As mentioned, the cases $q=-1$ and $q=0$
respectively correspond to the methods of steepest-descent and
minimal residues.) The extension to $P(\cdot)$ polynomial in $A$
is straightforward (including also linear combinations with
$A^{-1}$), using (\ref{gammak}).

The minimisation of $(P(A)g_{k+1},g_{k+1})$, or the calculation of
$\gamma_k$ in (\ref{gammak}), requires the calculations of terms
of the form $(A^n g,g)$, with $n=q$ or $n=q+1,q+2$. As shown
below, they are easily obtained from evaluations of $g(\cdot)$ at
different points. Notice that this construction implies that {\em
one} iteration of the algorithm will require {\em several}
evaluations of $g(\cdot)$. The construction proposed below is not
necessarily the most economical one, and evaluations of $f(\cdot)$
and $g(\cdot)$ at different points could be combined to provide
more efficient evaluations of terms $(A^n g,g)$. Our objective
here is simply to show that the family of algorithms considered in
the paper is not of purely theoretical interest, and that other
algorithms than the steepest-descent and minimal residues could
also be considered in practice.

Let $(A^n g,g)$ be the term to be evaluated, $n\geq 1$, with
$g=g(x)$ the gradient at the current point $x$. Define $x^{(0)}=x$
and
$$
x^{(i+1)}=x^{(i)}-\beta g(x^{(i)}) \,, \ i\geq 0 \,,
$$
with $\beta$ a fixed positive number (for instance, $\beta$ can be
taken equal to the value of $\gamma$ at previous iteration of the
algorithm). We obtain
$$
g^{(i)}=g(x^{(i)})=(I-\beta A)^i g \,.
$$
Define $P_i=(g,g^{(i)})=(g,(I-\beta A)^i g)$. In matrix notation,
$\Pb_n=\Qb_n G_n$, where
$$
\Pb_n=(P_0,P_1,\ldots,P_n)\TT\,, \
G_n=((g,g),(Ag,g),\ldots,(A^ng,g))\TT
$$
and the entries of the $(n+1)\times(n+1)$ matrix $\Qb_n$ are the
binomial coefficients,
$$
\Qb_n=\left( \begin{array}{ccccc} 1 \ & \ & \ & \ & \\
1 \ & \ -\beta \  &          &          & \\
1 \ & \ -2\beta \ & \ \beta^2 \ &          & \\
1 \ & \ -3\beta \ & \ 3\beta^2 \ & \ -\beta^3 \ & \ \ldots \\
\vdots \ & \ \vdots \ & \ \vdots \ & \ \vdots \ & \ \vdots
\end{array} \right) \,.
$$
The value of $(A^ng,g)$ is then directly obtained from
$G_n=\Qb_n^{-1}\Pb_n$. The entries of $\Pb_n$, defined by
$P_i=(g,g^{(i)})$, are also obtained more economically from
$$
P_{2j}=(g^{(j)},g^{(j)})\,, \ P_{2j+1}=(g^{(j+1)},g^{(j)})\,.
$$
Therefore, the evaluation of
$\gamma_k=(P(A)Ag_{k},g_{k})/(P(A)A^2g_{k},g_{k})$, with
$P(\cdot)$ a polynomial of degree $q$, requires $\lceil q/2
\rceil+2$ gradient evaluations (including the one at
$x^{(0)}=x_k$).

\subsection{Renormalisation}

We can rewrite the iteration (\ref{iteration}) as
$$
(x_{k+1} - x^*) = (x_k - x^*) - \gamma_k g_k \,,
$$
with $g_k=g(x_k)=A(x_k-x^*)$, so that
$$
g_{k+1} = g_k - \gamma_k A g_k = g_k -
\frac{(P(A)Ag_{k},g_{k})}{(P(A)A^2g_{k},g_{k})} \, Ag_k \,.
$$
Define the renormalised variable
\begin{equation}
z(x) =  \frac{Bg(x)}{(P(A)Ag(x),g(x))^{1/2}} \,, \label{zk}
\end{equation}
with $B=[P(A)A]^{1/2}$, the positive square-root of $P(A)A$, so
that $(z(x),z(x))=1$. Also define $z_k=z(x_k)$,
\begin{equation}\label{mui}
    \mu_j^k=(A^j z_k,z_k)\,, \ j\in\mathbb{Z}\,,
\end{equation}
so that $\mu_0^k=1$ for any $k$ and
$\gamma_k=\mu_0^k/\mu_1^k=1/\mu_1^k$. We obtain
\begin{eqnarray*}
    z_{k+1} &=& \frac{Bg_{k+1}}{(P(A)Ag_{k+1},g_{k+1})^{1/2}} =
    \frac{(I-\gamma_k A)B g_k}{((I-\gamma_k A)B g_k,(I-\gamma_k A)B
    g_k)^{1/2}} \\
 && = \frac{(I-\gamma_k A)z_k}{((I-\gamma_k A)z_k,(I-\gamma_k A)z_k)^{1/2}}
 =\frac{(I-\gamma_kA) z_k}{(1-2
\gamma_k\mu_1^k+\gamma_k^2\mu_2^k)^{1/2}} \,,
\end{eqnarray*}
that is,
\begin{equation}
z_{k+1} = \frac{(I-A/\mu_1^k) z_k}{(\mu_2^k/(\mu_1^k)^2-1)^{1/2}}
\,. \label{ziter}
\end{equation}
This gives the updating formula for the moments
\begin{equation}\label{muiter}
    \mu_j^{k+1} = (A^j z_{k+1},z_{k+1}) = \frac{\mu_j^k - 2
    \mu_{j+1}^k/\mu_1^k +
    \mu_{j+2}^k/(\mu_1^k)^2}{\mu_2^k/(\mu_1^k)^2-1} \,.
\end{equation}

In the special case where $\SH=\mathbb{R}^d$ we can assume that
$A$ is already diagonalised, with eigenvalues
$0<\ml_1\leq\ml_2\leq \cdots \leq\ml_d$. We can then consider
$[z_k]_i^2$, with $[z_k]_i$ the $i$-th component of $z_k$, as a
mass on the eigenvalue $\ml_i$, with $\sum_{i=1}^d [z_k]_i^2 =
\mu_0^k=1$. Define the discrete probability measure $\nu_k$
supported on $(\ml_1,\ldots,\ml_d)$ by $\nu_k(\ml_i)=[z_k]_i^2$,
so that its $j$-th moment is $\mu_j^k$, $j\in\mathbb{Z}$. We can
then interpret (\ref{ziter}) as a transformation $\nu_k \ra
\nu_{k+1}$. The asymptotic behaviour of the sequence $(z_k)$
generated by (\ref{ziter}) was studied in \cite{Akaike59}, see
also \cite{Forsythe68} and Chapter 7 of \cite{PWZl2000}. The main
result is that, assuming $0<\ml_1<\ml_2\leq \cdots
\leq\ml_{d-1}<\ml_d$, the sequence $(z_k)$ converges to a
two-dimensional plane, spanned by the eigenvectors $e_1$, $e_d$
associated with $\ml_1$ and $\ml_d$. The attraction property is
stated more precisely in Section~\ref{S:asymptoticbehaviour}, also
in the Hilbert space case. It is already important to notice that
although the results in the references above were obtained for the
steepest-descent algorithm, the renormalisation (\ref{zk}), which
depends on the chosen $P(\cdot)$, makes them applicable to any
algorithm in the family considered. Also, using the
renormalisation just defined we easily obtain (non asymptotic)
results on the monotonicity of the algorithm along its trajectory.

\subsection{Monotonicity of a rate of convergence}
\label{S:monotonicity}

Consider the function $(P(A) g_{k+1},g_{k+1})$ that $\gamma_k$
minimizes, and compute the rate of convergence $r_k$ of the
algorithm at iteration $k$, defined by
\begin{equation}\label{rk}
r_k=\frac{(P(A) g_{k+1},g_{k+1})}{(P(A)g_k,g_k)} \,.
\end{equation}
Other rates of convergence will be considered in
Section~\ref{S:rate} where they will be shown to be asymptotically
equivalent to $r_k$. Direct calculation gives
$r_k=\mbox{$1-1/L_k$}$, with
$$
L_k=\mu_1^k\mu_{-1}^k
$$
where the moments $\mu_i^k$ are defined by (\ref{mui}). Also, from
(\ref{muiter}), $L_k$ satisfies
$$
L_{k+1}-L_k= \frac{\mu_1^k}{D_k^2} \, \det \Mb_k
$$
with
$$
D_k=\mu_2^k - (\mu_1^k)^2
$$
and
\begin{equation}\label{Mk}
    \Mb_k=\left(
\begin{array}{lll}
\mu_{-1}^k  & \mu_0^k &\mu_1^k \\
\mu_0^k  & \mu_1^k &\mu_2^k \\
\mu_1^k  & \mu_2^k &\mu_3^k
\end{array}
\right) \,.
\end{equation}
The moment matrix $\Mb_k$ is positive semi-definite so that $\det
\Mb_k\geq 0$, and thus $L_{k+1}\geq L_k$, that is, both $L_k$ and
the rate $r_k$ are non-decreasing along the trajectory followed by
the algorithm. When $\SH=\mathbb{R}^2$ ($d=2$), $\det \Mb_k=0$ and
$r_k$ is constant. When $d>2$ or $\SH$ is a Hilbert space, the
rate is monotonically increasing for a typical $x_0$, indeed, for
almost all $z_0=z(x_0)$ with respect to the uniform measure on the
unit sphere when $\SH=\mathbb{R}^d$. Notice that if the rate is
constant over two iterations ($\det \Mb_k=0$), then the measure
$\nu_k$ is supported on two points only, and the iteration
(\ref{ziter}) for the masses shows that this situation will
continue: the rate will thus remain constant for all subsequent
iterations.

Note that $L_k$ and $D_k$ are bounded (since $\nu_k$ has a bounded
support), respectively by $L^*$ and $D^*$, with
$L^*=(M+m)^2/(4mM)$ and $D^*=(M-m)^2/4$, see Lemma \ref{L:D*} in
Appendix A3. Therefore, since $L_k$ is non-decreasing it converges
to some limit, and
\begin{equation}\label{detMk}
    \det \Mb_k = \frac{(L_{k+1}-L_k) D_k^2}{\mu_1^k} \leq
\frac{(L_{k+1}-L_k)(D^*)^2}{m} \ra 0\,, \ k\ra\infty\,.
\end{equation}

In addition to $L_k$ and $r_k$ another quantity also turns out to
be non-decreasing along the trajectory. Consider
\begin{equation}\label{def-Dk}
    \frac{(P(A)Ag_{k+1},g_{k+1})}{(P(A)A(x_{k+1}-x_k),(x_{k+1}-x_k))}
= \frac{(P(A)Ag_{k+1},g_{k+1})}{\gamma_k^2(P(A)Ag_k,g_k)} =
\mu_2^k-(\mu_1^k)^2 = D_k\,.
\end{equation}
Direct calculation using (\ref{ziter}) gives
$$
D_{k+1}-D_k = \frac{1}{D_k^2} \, \det \Nb_k
$$
with
$$
    \Nb_k=\left(
\begin{array}{lll}
\mu_{0}^k  & \mu_1^k &\mu_2^k \\
\mu_1^k  & \mu_2^k &\mu_3^k \\
\mu_2^k  & \mu_3^k &\mu_4^k
\end{array}
\right) \,.
$$
Again, $\Nb_k$ is positive semi-definite and $\det \Nb_k \geq 0$
so that $D_k$ is also non-decreasing. It converges to some limit
and $\det \Nb_k$ converges to zero for the same reasons as above.

Substitution of $P(A)$ for a particular algorithm shows which
quantities are monotonic. For the steepest-descent algorithm,
$P(A)=A^{-1}$, $(A^{-1}g_k,g_k)=2[f(x_k)-f(x^*)]$, and thus the
ratios $r_k=[f(x_{k+1})-f(x^*)]/[f(x_k)-f(x^*)]$ and
$D_k=(g_{k+1},g_{k+1})/((x_{k+1}-x_k),(x_{k+1}-x_k))$ are
monotonically non-decreasing. For the method of minimal residues,
$P(A)=I$, and the ratios $r_k=(g_{k+1},g_{k+1})/$ $(g_k,g_k)$ and
$D_k=(Ag_{k+1},g_{k+1})/(A(x_{k+1}-x_k),(x_{k+1}-x_k))$ are
monotonically non-decreasing.

The monotonicity and boundedness of $L_k$ and $D_k$ makes them
suitable for studying the asymptotic behaviour of the algorithm.
This is developed in the next section.

\section{Asymptotic behaviour of gradient algorithms}
\label{S:asymptoticbehaviour}

Consider the case $\SH=\mathbb{R}^d$, and assume that the minimal
and maximal eigenvalues of $A$, $\ml_1=m$, $\ml_d=M$, are simple.
The attraction property can be stated as follows. Choose
$z_0=z(x_0)$, the renormalised variable defined by (\ref{zk}) at
the initial point $x_0$, such that $(z_0,e_1)>0$, $(z_0,e_d)>0$,
with $e_1$ and $e_d$ the eigenvectors associated with $\ml_1$ and
$\ml_d$ respectively. Then
$$
z_{2k} \ra \sqrt{p} \, e_1 + \sqrt{1-p} \, e_d\,, \ \ z_{2k+1} \ra
\sqrt{1-p} \, e_1 - \sqrt{p} \, e_d \ \mbox{ when } \
k\ra\infty\,,
$$
where $p$ is some number in $(0,1)$, see Section~\ref{S:stability}
concerning the range of possible values for $p$. This property,
stated in a more general framework in Theorem
\ref{T:attraction-Rd} below, has important consequences for the
asymptotic rate of convergence of the algorithm, see
Section~\ref{S:rate}. The proof of the attraction property relies
on the convergence of successive transformations of the
probability measures $\nu_k$ defined by $[z_k]_i^2$. The
approaches used in \cite{Akaike59,Forsythe68} to study this
convergence do not apply when $\SH$ is infinite dimensional, and
we shall present a more general proof in
Section~\ref{S:theorem-tansform-pm}. It differs somewhat from the
one in Chapter 7 of \cite{PWZl2000}, in particular in the choice
of the monotonic sequence, $(L_k)$ instead of $(D_k)$.

The attraction theorem in $\mathbb{R}^d$ can be stated as follows.
We can assume that $A$ is diagonalised, and the probability
measure $\nu_k$ is then discrete and puts mass $[z_k]_i^2$ at the
eigenvalue $\ml_i$. Notice that the updating rule (\ref{ziter}) is
identical for $[z_k]_i$ and $[z_k]_j$ associated with
$\ml_i=\ml_j$, and the corresponding masses can thus be summed. We
can therefore assume that all eigenvalues are different when
studying the evolution of $\nu_k$, see Theorem
\ref{T:transform-Rd}.

\begin{theorem}
\label{T:attraction-Rd} Let $A$ be a $d\times d$ symmetric matrix,
positive definite, with minimum and maximum eigenvalues $m$ and
$M$ such that $0<m<M<\infty$ and apply a $P$-gradient algorithm,
see Definition \ref{D:P-gradient}, for the minimisation of $f(x)$
given by (\ref{fx}), initialized at $x_0$, with $z_0=z(x_0)$, see
(\ref{zk}). Assume that
\begin{equation}\label{condition-z0}
    E_1 z_0\neq 0 \ \mbox{ and } \ E_d z_0\neq 0\,,
\end{equation}
where $E_1$ and $E_d$ are the orthogonal projectors on the
eigenspaces respectively associated with $\ml_1=m$ and $\ml_d=M$.
Then the asymptotic behaviour of the renormalised gradient
$z_k=z(x_k)$ is such that
$$
z_{2k} = \sqrt{p} \, u_{2k} + \sqrt{1-p} \, v_{2k} \,, \ z_{2k+1}
= \sqrt{1-p} \, u_{2k+1} - \sqrt{p} \, v_{2k+1} \,,
$$
with $\|u_n\|=\|v_n\|=1$ $\forall n$, $\|A u_n-mu_n\| \ra 0$, $\|A
v_n-Mv_n\| \ra 0$ as $n\ra\infty$, and $p$, some number in
$(0,1)$, depending on $z_0$.
\end{theorem}

The proof is omitted since we prove later a more general property
valid for $\SH$ a Hilbert space. A more precise result is obtained
when the eigenvalues $\ml_1$ and $\ml_d$ are simple: the vector
$z_d$ converges to the two-dimensional plane defined by the
eigenvectors $e_1$ and $e_d$ associated with $\ml_1$ and $\ml_d$.

\begin{coro}
\label{C:attraction-Rd} Let $A$ be a positive-definite symmetric
matrix with ordered eigenvalues
$$
0 < m=\ml_1 < \ml_2 \leq \cdots \leq \ml_{d-1} < \ml_d=M
$$
and let $e_1$, $e_d$ be the eigenvectors associated with $\ml_1$
and $\ml_d$ respectively. Apply a $P$-gradient algorithm, see
Definition \ref{D:P-gradient}, for the minimisation of $f(x)$
given by (\ref{fx}), initialized at $x_0$ such that $z_0\TT
e_1\neq 0$ and $z_0\TT e_d\neq 0$, with $z_0=z(x_0)$, see
(\ref{zk}). Then the algorithm attracts to the plane $\Pi$ spanned
by $e_1$ and $e_d$ in the following sense:
$$
w\TT z_k \ra 0\,, \ k\ra \infty
$$
for any nonzero vector $w \in \Pi^\perp$. Moreover, the sequence
$(z_k)$ converges to a two-point cycle.
\end{coro}

This corollary is a straightforward consequence of Theorem
\ref{T:attraction-Rd}: when $\ml_1$ and $\ml_d$ are simple, with
associated eigenvectors $e_1$ and $e_d$, $u_n$ and $v_n$ then
respectively tend to $e_1$ and $e_d$. The result easily
generalizes to the case when (\ref{condition-z0}) is not
satisfied. The algorithm then attracts to a two-dimensional plane
defined by the eigenvectors $e_i$ and $e_j$ associated with the
smallest and largest eigenvalues such that $z_0\TT e_i\neq 0$ and
$z_0\TT e_j\neq 0$.

We state now the attraction theorem in the more general case where
$\SH$ is a Hilbert space. The proof is given in Appendix A1.

\begin{theorem}
\label{T:attraction-general} Let $A$ be a bounded real symmetric
operator in a Hilbert space $\SH$, positive, with bounds $m$ and
$M$, such that $0<m<M<\infty$ and apply a $P$-gradient algorithm,
see Definition \ref{D:P-gradient}, for the minimisation of $f(x)$
given by (\ref{fx}), initialized at $x_0$, with $z_0=z(x_0)$, see
(\ref{zk}). Assume that $z_0$ is such that for any $\me$,
$0<\me<(M-m)/2$,
\begin{equation}
(E_{m+\me}z_0,z_0) >0 \ \mbox{ and } \ (E_{M-\me}z_0,z_0) <1 \,,
\label{Cz0}
\end{equation}
with $(E_\ml)$ the spectral family of projections associated with
$A$. The asymptotic behaviour of the renormalised gradient
$z_k=z(x_k)$ is such that
\begin{equation}
z_{2k} = \sqrt{p} \, u_{2k} + \sqrt{1-p} \, v_{2k} \,, \ z_{2k+1}
= \sqrt{1-p} \, u_{2k+1} - \sqrt{p} \, v_{2k+1} \,, \label{zuv}
\end{equation}
with $\|u_n\|=\|v_n\|=1$ $\forall n$, $\|A u_n-mu_n\| \ra 0$, $\|A
v_n-Mv_n\| \ra 0$ as $n\ra\infty$, and $p$, some number in
$(0,1)$, depending on $z_0$.
\end{theorem}

\section{A property of successive transformations of a probability measure}
\label{S:theorem-tansform-pm}

The two properties established in this section form the
cornerstones of the proofs of the theorems of previous section. We
consider first the case of a discrete measure with finite support,
which in terms of convergence of a $P$-gradient algorithm
corresponds to the case $\SH=\mathbb{R}^d$. The proof is given in
Appendix A2.

\begin{theorem}
\label{T:transform-Rd} Let $\nu_0$ be a discrete probability
measure on $\{\ml_1,\ldots,\ml_d\}$ with
$$
0 < m=\ml_1 < \ml_2 < \cdots < \ml_{d-1} < \ml_d=M < \infty \,.
$$
Let $[z_k]_i^2$ denote the mass placed at $\ml_i$ by $\nu_k$, that
is, $\nu_k(\ml_i)=[z_k]_i^2$. Consider the transformation $T: \
\nu_k \ra \nu_{k+1}$ defined by
\begin{equation}
[z_{k+1}]_i = \frac{(1-\ml_i/\mu_1^k)
[z_k]_i}{(\mu_2^k/(\mu_1^k)^2-1)^{1/2}} \label{nuk+1-Rd}
\end{equation}
with the moments $\mu_i^k$ defined by (\ref{mui}). Then, when
$k\ra\infty$,
\begin{equation}
([z_{2k}]_1)^2 \ra p \,, \ ([z_{2k+1}]_1)^2 \ra 1-p \ \mbox{ and }
([z_{2k}]_d)^2 \ra 1-p \,, \ ([z_{2k+1}]_1)^2 \ra p
\label{nukasymptotic-Rd}
\end{equation}
for some $p$ depending on $\nu_0$, $0<p<1$. Furthermore,
$$
p= \frac{1}{2} \pm \frac{\rho+1}{\rho-1} \sqrt{\frac{1}{4}-
\frac{\rho L}{(\rho+1)^2}}
$$
with $\rho=M/m$ and $L=\lim_{k\ra\infty} \mu_1^k\mu_{-1}^k$.
\end{theorem}

Note that the limiting value $L$ depends on $\nu_0$, so that the
value of $p$ that characterizes the attractor is difficult to
predict. The range of possible values for $p$ is discussed in
Section~\ref{S:stability}.

We consider now the case of an arbitrary measure on an interval,
which raises some additional difficulties compared to previous
case. In terms of convergence of a $P$-gradient algorithm, it
corresponds to the case where $\SH$ is a Hilbert space: for
$E_\ml$ the spectral family associated with the operator $A$, we
define the measure $\nu_k$ by $\nu_k(d\ml)=d(E_\ml z_k,z_k)$,
$m\leq\ml\leq M$.

\begin{theorem}
\label{T:transform} Let $\nu_0$ be a probability  measure on the
family $\SB$ of Borel sets of $(0,\infty)$, with support $[m,M]$,
so that
\begin{eqnarray*}
m &=& \essinf(\nu_0) = \sup(\alpha \ / \ \nu_0\{x\,, \, x<\alpha\}=0) \,, \\
M &=& \esssup(\nu_0) = \inf(\alpha \ / \ \nu_0\{x\,, \,
x>\alpha\}=0) \,.
\end{eqnarray*}
Assume that $0<m<M<\infty$. Consider the transformation $T: \
\nu_k \ra \nu_{k+1}$ defined by
\begin{equation}
\nu_{k+1}(\SA) = \int_\SA \frac{(\ml-\mu_1^k)^2}{D_k} \,
\nu_k(d\ml) \label{nuk+1}
\end{equation}
for any $\SA\in\SB$, where $\mu_1^k=\int\ml \, \nu_k(d\ml)$ and
$D_k=\mu_2^k-(\mu_1^k)^2$, with $\mu_2^k=\int \ml^2 \,
\nu_k(d\ml)$. Then, as $k\ra\infty$,
\begin{equation}
\nu_{2k}(\SI) \ra p \,, \ \nu_{2k+1}(\SI) \ra 1-p
\label{nukasymptotic}
\end{equation}
for all $\SI=[m,x)$, $m<x<M$, for some $p$ depending on $\nu_0$,
$0<p<1$.
\end{theorem}

The proof of Theorem \ref{T:transform} is given in Appendix A3.
\section{Stability of attractors}
\label{S:stability}

The range of possible values for $p$ in the attraction Theorem
\ref{T:attraction-Rd} ($\SH=\mathbb{R}^d$) is considered in
Theorem 3 of \cite{Akaike59} (see also Lemma 3.5 of
\cite{NocedalSZ2002}). Let $s(\ml)$ and $\ml^*$ be defined by
(\ref{sl}). This theorem states that when $\ml^*$ is not discarded
at any iteration, that is, when $\mu_1^k \neq \ml^*$ for any $k$,
then $p\in[1/2-s(\ml^*),1/2+s(\ml^*)]$ (note that this assumption
cannot be checked). In this section we extend this result in two
directions: (i) we will assume that $\SH$ is a Hilbert space, (ii)
we study the stability of the attractor defined by $p$ in Theorem
\ref{T:attraction-general}. We shall use the following definition
of stability, see \cite{HaleK91} p.\ 444, \cite{LaSalle76}, p.\ 7.

\begin{definition}
A fixed point $\nu^*$ for a mapping $T(\cdot)$ on a metric space
with distance $d(\cdot,\cdot)$ will be called stable if $\forall
\me>0$, $\exists \alpha>0$ such that for any $\nu_0$ for which
$d(\nu_0,\nu^*) < \alpha$, $d(T^n(\nu_0),\nu^*)<\me$ for all
$n>0$. A fixed point $\nu^*$ is unstable if it is not stable.
\end{definition}

We shall use the distance $d(\nu,\nu')$ given by the
L\'evy-Prokhorov metric, see \cite{Shiryaev96} p.\ 349. In our
case (measures supported on $[m,M]$), $d(\nu,\nu')$ becomes the
L\'evy distance between the distribution functions $F,F'$
associated with $\nu,\nu'$, which we denote
$$
L(F,F')=\inf \{\me: F'(x-\me)-\me \leq F(x) \leq F'(x+\me)+\me\,,
\; \forall x \} \,.
$$
In the case where one of the two measures is the discrete measure
$\nu_p^*$  concentrated on $m,M$, with $\nu_p^*(m)=p$,
$\nu_p^*(M)=1-p$, we get
\begin{eqnarray*}
d(\nu,\nu_p^*) &=& L(F,F_p^*) \\
&& \hspace{-1cm} = \inf\{\me: F(x)\leq p+\me \mbox{ for } x<M-\me
\; \mbox{ and } \; p-\me \leq F(x) \mbox{ for } m+\me \leq x \}
\,,
\end{eqnarray*}
with $F_p^*$ the distribution function associated with $\nu_p^*$.
We then have proved the following, see Appendix A4.

\begin{theorem}
\label{T:instability} Consider the situation of Theorem
\ref{T:transform}, with $\nu_0$ any probability measure supported
on some closed subset
 $\SS_A$ of $[m,M]$ and
 $$
 \essinf(\nu_0)=m\,, \ \esssup(\nu_0)=M\,.
 $$
(i) The measure $\nu_p^*$ is a fixed point for the mapping $T^2$.\\
(ii) Consider the set $\SI_u$ defined by
$$
\SI_u=\left(0,\frac{1}{2}-s(\ml^*)\right) \cup
\left(\frac{1}{2}+s(\ml^*),1\right) \,,
$$
where
\begin{equation}
s(\ml)=\frac{\sqrt{(M-\ml)^2+(\ml-m)^2}}{2(M-m)} \,, \
\ml^*=\min_{\ml\in\SS_{A}} s(\ml) \,. \label{sl}
\end{equation}
Any fixed point $\nu_p^*$ with $p$ in $\SI_u$ corresponds to an unstable fixed point for $T^2$.\\
(iii) Any point in the interval
\begin{equation}\label{Is}
\SI_s=\left(\frac{1}{2}-s(\ml^*), \frac{1}{2}+s(\ml^*)\right)
\end{equation}
corresponds to a stable $\nu_p^*$ for the mapping $T^2$.
\end{theorem}

\begin{remark}
The convergence $d(\nu_k,\nu_p^*)\ra 0$ is equivalent to weak
convergence $\nu_k \stackrel{\rm w}{\longrightarrow} \nu_p^*$ in
the usual sense. If $z_k$ is associated with the spectral measure
$\nu_k$ and $z_p^*$ with $\nu_p^*$, then, in the Hilbert space
this is equivalent to $(z_k-z_p^*,y) \ra 0$ for any $y\in\SH$,
whereas strong convergence would require $\|z_k-z_p^*\| \ra 0$.
For $\mathbb{R}^d$, the two types of convergence are equivalent,
and thus Corollary \ref{C:attraction-Rd} implies strong
convergence. However, for $\SH$ a Hilbert space the equivalence is
false, and indeed strong convergence generally does {\em not}
hold. The stability property (iii) is thus a weak statement when
$\SH$ is a Hilbert space. The $L_2$ metric in $\SH$ induces the
Hellinger metric on the space of spectral measures, which defines
the same topology as the distance in variation, see
\cite{Shiryaev96}, p.\ 364. Strong convergence in $\SH$ is thus
related to distance in variation in the space of spectral measures
and is clearly difficult to obtain --- except in the special
situation where $\nu_0$ has positive mass at $\{m\}$ and $\{M\}$
and presents a spectral gap: $\nu_0[(m,m+\me)]=0$ and
$\nu_0[(M-\me,M)]=0$ for some $\me>0$.
\end{remark}

We have $\nu_{k+2}(d\ml)=H(\nu_k,\ml)\nu_k(d\ml)$, with
$H(\nu_k,\ml)$ given by (\ref{Hnuk}) in Appendix A4. One may then
notice that when $\nu_0$ is a discrete probability measure, the
condition $H(\nu_p^*,\ml)>1$ used in the proof of the instability
part of the theorem, see Appendix A4, corresponds to a condition
on the eigenvalues of the Jacobian of the transformation $T^2$,
see \cite{PWZl2000}.

\begin{figure}
\centering
\resizebox{0.75\textwidth}{!}{%
  \includegraphics{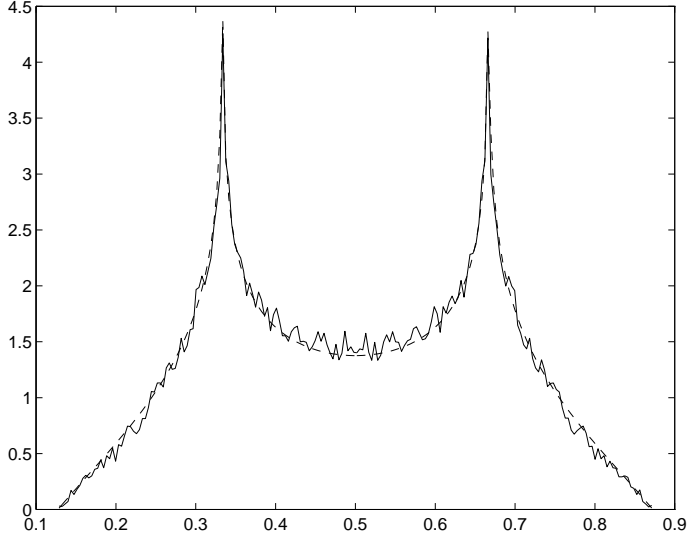}
}
\caption{Empirical density of attractors (full line) and
$\varphi(p)$, see (\ref{density}), for $d=3$ ($m=1$, $\ml=4$,
$M=10$)} \label{F:density}
\end{figure}

Note that the stability interval $\SI_s$ always contains the
interval
$$
\left(\frac{1}{2}-\frac{1}{2\sqrt{2}},\
\frac{1}{2}+\frac{1}{2\sqrt{2}}\right) \ \approx \ (0.14645, \
0.85355) \,.
$$
Numerical simulations for $\SH=\mathbb{R}^3$, with $A$ having
eigenvalues $m<\ml<M$, show that for any initial density of $x_0$
in $\mathbb{R}^d$ associated with a density of $z_0$ reasonably
spread on the unit sphere, the density of the values of $p$
corresponding to stable attractors $\nu_p^*$ can be approximated
by
\begin{equation}\label{density}
\varphi(p)=C \log [\min\{1,H(\nu_p^*,\ml)\}] = \left\{
\begin{array}{cl}
C \log H(\nu_p^*,\ml) & \mbox{ if } p \in \SI_s \\
0 & \mbox{ otherwise}\,,
\end{array} \right.
\end{equation}
where $C$ is a normalisation constant and $H(\nu_p^*,\ml)$ is
given by (\ref{Hnup*}). Figure~\ref{F:density} shows the empirical
density of attractors (full line) together with $\varphi(p)$
(dashed line) in the case $m=1$, $\ml=4$, $M=10$. The support of
this density coincides with the stability interval $\SI_s$ given
by (\ref{Is}). When $d>3$, the density of attractors depends on
the initial density of $x_0$.

\section{Rates of convergence}
\label{S:rate}

We first state a property showing that different definitions of
rates of convergence are asymptotically equivalent, see Appendix
A5 for the proof.

\begin{theorem}
\label{T:rates-equivalent} Let $W$ be a bounded positive
self-adjoint operator in $\SH$, with bounds $c$ and $C$ such that
$0<c<C<\infty$. Assume that $W$ commutes with $A$ (when
$\SH=\mathbb{R}^d$, $W$ is a $d\times d$ positive-definite matrix
with minimum and maximum eigenvalues respectively $c$ and $C$).
Define
$$
R_k(W) = \frac{(Wg_{k+1},g_{k+1})}{(Wg_k,g_k)}
$$
if $\|g_k\|\neq 0$ and $R_k(W)=1$ otherwise. Apply a $P$-gradient
algorithm (\ref{iteration}), initialized at $x_0$, with $\gamma_k$
given by (\ref{gammak}), for the minimisation of $f(x)$ given by
(\ref{fx}), with minimum value at $x^*$. Then the limit
$$
R(W,x_0,x^*) = \lim_{n\ra\infty} \left[ \prod_{k=0}^{n-1} R_k(W)
\right]^{1/n}
$$
exists for all $x_0,x^*$ in $\SH$ and $R(W,x_0,x^*)=R(x_0,x^*)$
does not depend on $W$. In particular,
$$
R(W,x_0,x^*)=\lim_{n\ra\infty} \left( \prod_{k=0}^{n-1} r_k
\right)^{1/n}
$$
with $r_k$ defined by (\ref{rk}).
\end{theorem}

\begin{figure}
\centering
\resizebox{0.75\textwidth}{!}{%
  \includegraphics{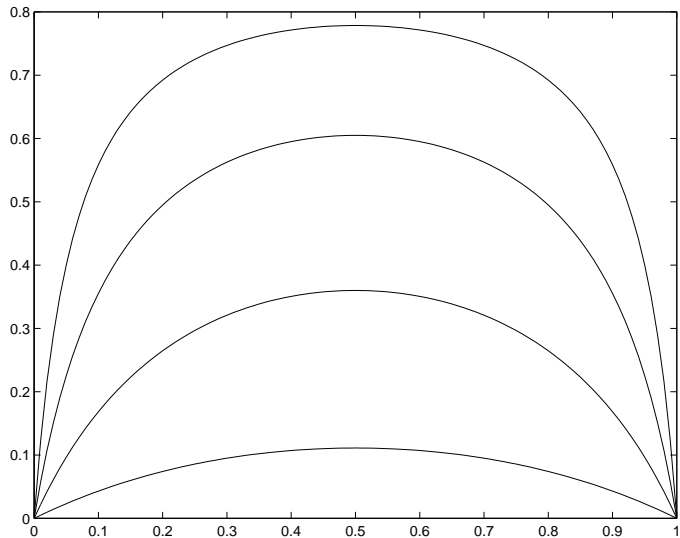}
}
\caption{$r(p)$ as a function of $p$, for $\rho=2$ (bottom
curve), 4, 8  and 16 (top)} \label{F:rates-p}
\end{figure}

From the results of Section~\ref{S:asymptoticbehaviour}, we have
$$
R(W,x_0,x^*) = r(p) =
\frac{p(1-p)(\rho-1)^2}{[p+\rho(1-p)][(1-p)+\rho p]}
$$
for any $W$, where $p$ defines the attractor, see (\ref{zuv}), and
$\rho=M/m$ is the condition number of the operator. The function
$r(p)$ is symmetric with respect to $1/2$ and monotonously
increasing from 0 to $1/2$, see Figure~\ref{F:rates-p}. The worst
asymptotic rate is thus obtained at $p=1/2$:
\begin{equation}
R_{\max}=\left( \frac{\rho-1}{\rho+1} \right)^2 \,. \label{Rmax}
\end{equation}

\begin{figure}
\centering
\resizebox{0.75\textwidth}{!}{%
  \includegraphics{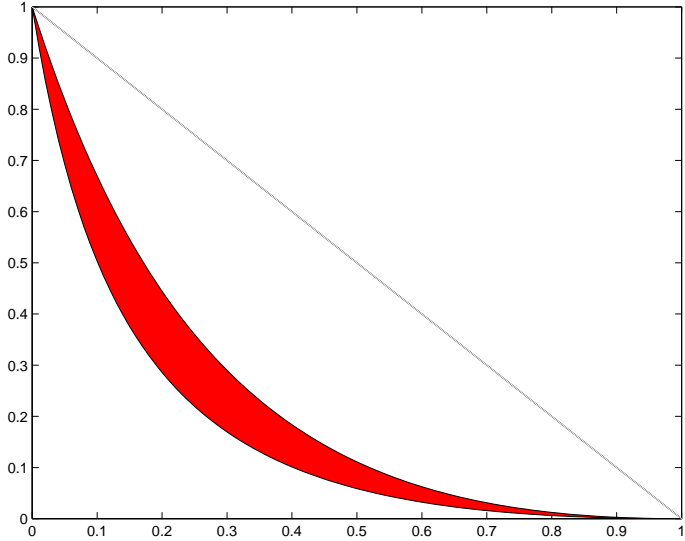}
}
\caption{Range $[R_{\min}^*,R_{\max}]$ of possible values of the
asymptotic rate $r(p)$ as a function of $1/\rho$} \label{F:rates}
\end{figure}

Note that $\forall k$,  $r_k \leq R_{\max}$ since $r_k$ is not
decreasing, see Section~\ref{S:monotonicity}. For a typical $x_0$
(such that the convergence is not finite, that is, such that
$r(p)\neq 0$), the stability analysis of Section~\ref{S:stability}
shows that only values of $p$ in $\SI_s$ given by (\ref{Is}) may
correspond to stable attractors. The range of possible values of
$R(p)$ is thus $[R_{\min},R_{\max}]$, where $R_{\max}$, given by
(\ref{Rmax}), is obtained for $p=1/2$ and
$$
R_{\min} \leq R_{\min}^* = R(1/2+1/[2\sqrt{2}]) =
\frac{(\rho-1)^2}{(\rho+1)^2+4\rho} \,.
$$
Figure~\ref{F:rates} presents the range $[R_{\min}^*,R_{\max}]$ as
a function of $1/\rho$, the upper curve corresponding to
$R_{\max}$ and the lower to $R_{\min}^*$. The maximum size of the
range is $3-2\sqrt{2}\simeq 0.1716$, obtained at
$\rho=1+2\sqrt{2}+2\sqrt{2+\sqrt{2}} \simeq 7.5239$. These results
confirm the experimental observation that the rate of convergence
of the gradient algorithm is generally close to its worst value
$R_{\max}$, see \cite{NocedalSZ2002}. The same property is true
for any $P$-gradient algorithm.

\begin{remark}
A similar analysis for $D_k$ defined by (\ref{def-Dk}), which is
also not decreasing, shows that $D_k \ra D(p)=p(1-p)(M-m)^2$ as
$k\ra\infty$, with $D_k\leq D^*=D(1/2)=(M-m)^2/4$ for all $k$.
Also, for any typical $x_0$ such that $p\in\SI_s$ given by
(\ref{Is}), we have $D(p)\geq D(1/2+1/[2\sqrt{2}])=(M-m)^2/8$.
\end{remark}

Another quantity of interest is given by
$$
\Delta_N = \log(R_{\max}/R_{\min})/[\log(R_{\max})\log(R_{\min})]
\,.
$$
Indeed, for $N$ large enough, $(Wg_N,g_N)/(Wg_0,g_0) \simeq
r(p)^N$, the number $N$ of iterations required for obtaining a
ratio $(Wg_N,g_N)/(Wg_0,g_0)=\me$ $(\me\ll 1)$ is approximately
$\log(\me)/\log[r(p)]$ and $\Delta_N |\log(\me)|$ thus indicates
the length of the interval of possible values for $N$ due to the
range of possible values for $p$. Direct calculation gives
$\Delta_N\, |\log(R_{\max})|<1/2$ for any $\rho$ and
$$
\Delta_N = \rho/8-1/4+\SO(1/\rho)\,, \ 1/\log(R_{\max}) =
-\rho/4+\SO(1/\rho)
$$
for large $\rho$. Therefore, the number of iterations required by
a $P$-gradient algorithm to achieve a given precision $\me<<1$
varies at most by a factor 2 depending on the (typical) starting
point $x_0$, factors of variation close to 2 being possible only
when $\rho$ is large.

The average value of $R(W,x_0,x^*)$ for $z_0=z(x_0)$ uniformly
distributed on the unit sphere is the same for any $P$-gradient
algorithm, more generally, the distribution of $R(W,x_0,x^*)$
associated with a particular distribution of $z_0$ does not depend
on the particular $P$-gradient algorithm considered. Moreover,
numerical simulations show that the average value of
$R(I,x_0,x^*)$ is the same for the steepest-descent
($P(A)=A^{-1}$) and minimum residues ($P(A)=I$) algorithms for
$x_0$ uniformly distributed on the sphere $\|x_0-x^*\|=1$. The
small deviations in average performance between different
$P$-gradient algorithms can only be related to the fact that a
fixed distribution for $x_0$ corresponds to different
distributions for $z(x_0)$.

\begin{remark}
It is known that the introduction of a relaxation coefficient
$\gamma$, with $0<\gamma<1$, in the steepest-descent algorithm
totally changes its behaviour, see, e.g., Chapter 7 of
\cite{PWZl2000}; the algorithm (\ref{algo}) then becomes $x_{k+1}
= x_k - \gamma [(g_k,g_k)/(Ag_k,g_k)] g_k$. For $\SH=\mathbb{R}^d$
and a fixed $A$, depending on the value of $\gamma$, the
renormalized process either converges to periodic orbits (the same
for almost all starting points) or exhibits a chaotic behaviour,
with the classical period-doubling phenomenon in the case $d=2$.
In higher dimensions, repeated numerical trials show that the
process typically no longer converges to the 2-dimensional plane
spanned by $(e_1,e_d)$. A detailed analysis for $d=2$ and
experimental results for $d>2$ also show that relaxation (with
$\gamma$ close to 1) considerably improves the rate of
convergence. Similar results hold more generally for all
$P$-gradient algorithms, with the iteration (\ref{iteration})
transformed into $x_{k+1} = x_k -\gamma\gamma_k g_k$, with
$\gamma$ the (fixed) relaxation coefficient and $\gamma_k$ given
by (\ref{gammak}). Steepest descent with random relaxation
coefficient $\mg\in(0,2)$ is considered in \cite{RaydanS2002},
avoiding the two point attraction and significantly improving the
behavior of ordinary steepest descent.
\end{remark}

\section*{Appendix}

\noindent {\bf A1.} {\em Proof of Theorem
\ref{T:attraction-general}.} The proof relies on Theorem
\ref{T:transform} (Theorem \ref{T:transform-Rd} when
$\SH=\mathbb{R}^d$), which concerns successive transformations
applied to a probability measure.

Since $A$ is self-adjoint, its spectrum $\SS_A$ is a closed subset
of the interval $[m,M]$ of the real line and $m,M\in\SS_A$. Let
$E_\ml$ be the spectral family associated with $A$, and define the
spectral measure $\nu_k$ by $\nu_k(d\ml)=d(E_\ml z_k,z_k)$,
$m\leq\ml\leq M$. Since $(z_k,z_k)=\int_m^M \nu_k(d\ml)=1$,
$\nu_k$ is a probability measure on the Borel sets of
$(0,\infty)$, with $\nu_k([m,M])=1$ $\forall k$. This
representation gives
$$
\mu_1=(Az_k,z_k)=\int \ml \, \nu_k(d\ml)\,, \
\mu_2=(A^2z_k,z_k)=\int \ml^2 \, \nu_k(d\ml)
$$
where integration is over $[m,M]$ unless otherwise specified.
Therefore, for any Borel set $\SA$ the transformation
(\ref{ziter}) gives in terms of $\nu_k$:
$$
\nu_{k+1}(\SA) = \frac{\int_\SA \left[\ml-\int \ml' \,
\nu_k(d\ml')\right]^2 \, \nu_k(d\ml)} {\int \ml'^2 \, \nu_k(d\ml')
- \left[ \int \ml' \, \nu_k(d\ml') \right]^2 }  \,.
$$
The conditions (\ref{Cz0}) on $z_0$ are equivalent to
$\essinf(\nu_0)=m$ and $\esssup(\nu_0)=M$, see Theorem
\ref{T:transform}, and the updating rule for $\nu_k$ can be
written as (\ref{nuk+1}). Theorem \ref{T:transform} then implies
(\ref{nukasymptotic}), which can be written as: $\forall \me>0$,
$\me\leq \beta=(M-m)/2$,
\begin{eqnarray*}
(E_{m+\me}z_{2k},z_{2k}) \ra p \,, && (E_{M-\me}z_{2k},z_{2k}) \ra p \,, \\
(E_{m+\me}z_{2k+1},z_{2k+1}) \ra 1-p \,, &&
(E_{M-\me}z_{2k+1},z_{2k+1}) \ra 1-p \,,
\end{eqnarray*}
as $k\ra\infty$, where $p$ depends on $z_0$, $0<p<1$. Define
$p_{2k}=(E_{m+\beta}z_{2k},z_{2k})$,
$p_{2k+1}=1-(E_{m+\beta}z_{2k+1},z_{2k+1})$, and the angles
$\varphi$, $\varphi_n$ by $\cos \varphi=\sqrt{p}$, $\sin
\varphi=\sqrt{1-p}$, $\cos \varphi_n=\sqrt{p_n}$, $\sin
\varphi_n=\sqrt{1-p_n}$, $\forall n$. Also define $s_{2k}=
E_{m+\beta}z_{2k}/\cos \varphi_{2k}$, $s_{2k+1}=
E_{m+\beta}z_{2k+1}/\sin \varphi_{2k+1}$, $t_{2k}=
(z_{2k}-E_{m+\beta}z_{2k})/\sin \varphi_{2k}$, $t_{2k+1}=
-(z_{2k+1}-E_{m+\beta}z_{2k+1})/\cos \varphi_{2k+1}$. This gives
$p_n \ra p$ as $n\ra\infty$, $\|s_n\|=\|t_n\|=1$ $\forall n$, and
$z_{2k}=\cos \varphi_{2k} \, s_{2k} + \sin \varphi_{2k} \,
t_{2k}$, $z_{2k+1}=\sin \varphi_{2k+1} \, s_{2k+1} - \cos
\varphi_{2k+1} \, t_{2k+1}$. Also,
$$
\|As_n-ms_n\|^2 = \int (\ml-m)^2 \, d(E_\ml s_n,s_n)\,,
$$
which, for $n=2k$ and any $\me$, $0<\me<\beta$, gives
\begin{eqnarray*}
\|As_{2k}-ms_{2k}\|^2 &=& \int_m^{m+\beta} \frac{(\ml-m)^2}{p_{2k}} \, d(E_\ml z_{2k},z_{2k}) \\
&=& \int_m^{m+\me} \frac{(\ml-m)^2}{p_{2k}} \, d(E_\ml
z_{2k},z_{2k}) +
\int_{m+\me}^{m+\beta} \frac{(\ml-m)^2}{p_{2k}} \, d(E_\ml z_{2k},z_{2k}) \\
&&\leq \frac{\me^2}{p_{2k}} + \frac{\beta^2}{p_{2k}}
\left[p_{2k}-\int_m^{m+\me} d(E_\ml z_{2k},z_{2k})\right] \,.
\end{eqnarray*}
Since $p_{2k}\ra p$ and $\int_m^{m+\me} d(E_\ml z_{2k},z_{2k}) \ra
p$ as $k\ra\infty$, $\|As_{2k}-ms_{2k}\| \ra 0$ as $k\ra\infty$.
Similarly, $\|As_{2k+1}-ms_{2k+1}\| \ra 0$ as $k\ra\infty$ and
$\|At_{n}-Mt_{n}\| \ra 0$ as $n\ra\infty$. Consider now
$$
u_n=\cos \vartheta_n \, s_n + \sin \vartheta_n \, t_n \,, \ v_n=
-\sin \vartheta_n \, s_n + \cos \vartheta_n \, t_n \,.
$$
Straightforward calculations show that
$\vartheta_n=\varphi_n-\varphi$ gives (\ref{zuv}) with
$\|u_n\|=\|v_n\|=1$ $\forall n$. Also
$$
\|Au_n-mu_n\| \leq |\cos \vartheta_n| \|As_n-ms_n\| + |\sin
\vartheta_n| (M-m) \,,
$$
and, since $\|As_n-ms_n\| \ra0$, $\vartheta_n \ra 0$ as
$n\ra\infty$, $\|Au_n-mu_n\| \ra 0$ as $n\ra\infty$. Similarly,
$\|Av_n-Mv_n\| \ra0$ as $n\ra\infty$. \carre

\vspace{1cm} \noindent {\bf A2.} {\em Proof of Theorem
\ref{T:transform-Rd}.} We first prove that the mass of $\nu_k$
tends to concentrate on two eigenvalues only. When $\nu_0$ is non
degenerate, $L_1>1$ from Jensen inequality, and thus, since
$(L_k)$ is non-decreasing, see Section~\ref{S:monotonicity}, $L_k
\geq L_1 >1$. Now, from Lagrange identity $(\sum a_i^2)(\sum
b_i^2) = \sum_{i<j}(a_ib_j-a_jb_i)^2+(\sum a_ib_i)^2$
\begin{eqnarray*}
  L_k &=& \left( \sum_{i=1}^d \ml_i [z_k]_i^2 \right) \, \left(
\sum_{i=1}^d [z_k]_i^2/\ml_i \right) \\
   &=& \sum_{i<j} [z_k]_i^2 [z_k]_j^2
   \left( \frac{\sqrt{\ml_i}}{\sqrt{\ml_j}} - \frac{\sqrt{\ml_j}}{\sqrt{\ml_i}} \right)^2 +
   \left(\sum_{i=1}^d [z_k]_i^2 \right)^2 \\
   &=& \sum_{i<j} [z_k]_i^2 [z_k]_j^2
   \frac{(\ml_i-\ml_j)^2}{\ml_i\ml_j} + 1  \,.
\end{eqnarray*}
Let $i_k$ and $j_k$ denote the indices that achieve $\max_{i<j}
[z_k]_i^2 [z_k]_j^2$. We have
$$
L_k \leq [z_k]_{i_k}^2 [z_k]_{j_k}^2 \sum_{i<j}
\frac{(\ml_i-\ml_j)^2}{\ml_i\ml_j} + 1
$$
and thus
$$
[z_k]_{i_k}^2 [z_k]_{j_k}^2 \geq \delta = \frac{L_1-1}{\sum_{i<j}
\frac{(\ml_i-\ml_j)^2}{\ml_i\ml_j}} \,.
$$
Moreover, $[z_k]_{i_k}^2+[z_k]_{j_k}^2< 1$ gives
$$
\delta <[z_k]_{i_k}^2 < 1-\delta \ \mbox{ and } \ \delta
<[z_k]_{j_k}^2 < 1-\delta \,.
$$
Consider the matrix $\Mb_k$ given by (\ref{Mk}). Its determinant
can be written as
\begin{eqnarray*}
\det \Mb_k &=& \sum_{i<j<l} [z_k]_i^2 [z_k]_j^2 [z_k]_l^2
\frac{(\ml_i-\ml_j)^2(\ml_i-\ml_l)^2(\ml_j-\ml_l)^2}{\ml_i\ml_j\ml_l}
\\
&& \geq [z_k]_{i_k}^2 [z_k]_{j_k}^2 (\ml_{i_k}-\ml_{j_k})^2
\sum_{i \neq i_k, j_k} [z_k]_i^2
\frac{(\ml_i-\ml_{i_k})^2(\ml_i-\ml_{j_k})^2}{\ml_i\ml_{i_k}\ml_{j_k}}
\\
&& \geq \delta \frac{\delta_\ml^6}{M^3} \sum_{i \neq i_k, j_k}
[z_k]_i^2
\end{eqnarray*}
where
$$
\delta_\ml = \min_{i,j} |\ml_i-\ml_j| \,.
$$
Since $\det \Mb_k \ra 0$ as $k\ra\infty$, see (\ref{detMk}), we
get $\sum_{i \neq i_k, j_k} [z_k]_i^2 \ra 0$ as $k\ra\infty$. The
mass thus tends to concentrate on $\ml_{i_k}, \ml_{j_k}$.

Next we prove that $i_k$ and $j_k$ eventually become fixed. From
the result above, $\forall \me>0$, $\exists k_\me$ such that
$\sum_{i \neq i_k, j_k} [z_k]_i^2 <\me$, $k> k_\me$.

Consider the updating equation (\ref{nuk+1-Rd}). We have for any
$i$, $(\mu_1^k-\ml_i)^2 \leq (M-m)^2$. Also,
$D_k=\mu_2^k-(\mu_1^k)^2 \geq D_0$, see
Section~\ref{S:monotonicity}. This gives for $i\neq i_k,j_k$ and
$k>k_\me$
$$
[z_{k+1}]_i^2 < \me \frac{(M-m)^2}{D_0} \,.
$$
Taking $\me_1=\delta D_0/(M-m)^2$ we obtain $[z_{k+1}]_i^2 <
\delta$ for $i\neq i_k,j_k$ and $k>k_{\me_1}$. Since
$[z_{k+1}]_{i_{k+1}}^2>\delta$ and $[z_{k+1}]_{j_{k+1}}^2>\delta$,
$i\not\in \{i_k,j_k\}$ implies $i\not\in \{i_{k+1},j_{k+1}\}$,
$k>k_{\me_1}$ and thus $\{i_k,j_k\}=\{i^*,j^*\}$ for
$k>k_{\me_1}$.

We show now that $\{i^*,j^*\}=\{1,d\}$. Assume that $i^*<j^*<d$
(which implies $[z_k]_d^2 \ra 0$, $k\ra\infty$). We need to show
that $(\ml_d-\mu_1^k)^2>(\ml_{j^*}-\mu_1^k)^2$ for $k$ large
enough. We have
$$
\mu_1^k = \ml_{i^*}[z_k]_{i^*}^2 + \ml_{j^*}[z_k]_{j^*}^2 +
\sum_{i\neq i^*,j^*} \ml_i [z_k]_i^2 \leq \ml_{i^*}[z_k]_{i^*}^2 +
\ml_{j^*}[z_k]_{j^*}^2 + \ml_d \sum_{i\neq i^*,j^*} [z_k]_i^2 \,.
$$
Take $\me_2=\min\{\me_1, \delta \delta_\ml/\ml_d\}$. For
$k>k_{\me_2}$ we have
$$
\mu_1^k \leq \ml_{i^*}[z_k]_{i^*}^2 + \ml_{j^*}[z_k]_{j^*}^2 +
\ml_d \me_2 \leq \ml_{i^*} \delta + \ml_{j^*} (1-\delta) + \ml_d
\me_2 \leq \ml_{j^*} - \delta\delta_\ml + \ml_d \me_2 \leq
\ml_{j^*}
$$
and thus $(\ml_d-\mu_1^k)^2>(\ml_{j^*}-\mu_1^k)^2$. From
(\ref{nuk+1-Rd}), this gives for $k>k_{\me_2}$
$$
\left(\frac{[z_{k+1}]_d}{[z_k]_d}\right)^2 = \frac{
(\ml_d-\mu_1^k)^2 }{D_k}
> \frac{ (\ml_{j^*}-\mu_1^k)^2 }{D_k} =
\left(\frac{[z_{k+1}]_{j^*}}{[z_{k}]_{j^*}}\right)^2
$$
and thus
$$
\left(\frac{[z_{k+1}]_{j^*}}{[z_{k+1}]_{d}}\right)^2 <
\left(\frac{[z_{k}]_{j^*}}{[z_{k}]_{d}}\right)^2 \,.
$$
We arrived at a contradiction since $[z_k]_d^2 \ra 0$ and
$[z_{k}]_{j^*}^2$ is bounded from below by $\delta$. Therefore
$j^*=d$. Similarly, $i^*=1$.

Finally, let $L$ denote $\lim_{k\ra\infty} L_k$, see
Section~\ref{S:monotonicity}. There are only two discrete measures
with nonzero masses on $\ml_1$ and $\ml_d$ and such that
$\mu_1\mu_{-1}=L$,
$$
\nu^{(1)} = \left\{ \begin{array}{cc} \ml_1 & \ml_d \\ p & 1-p
\end{array} \right\} \ \mbox{ and } \ \nu^{(2)} = \left\{ \begin{array}{cc} \ml_1 & \ml_d \\ 1-p & p
\end{array} \right\}
$$
with
$$
p= \frac{1}{2} - \frac{\rho+1}{\rho-1} \sqrt{\frac{1}{4}-
\frac{\rho L}{(\rho+1)^2}}
$$
and $\rho=M/m$. Direct calculation shows that $\nu_k=\nu^{(1)}$
gives $\nu_{k+1}=\nu^{(2)}$, hence the convergence of $\nu_k$ to
the cyclic attractor $\nu^{(1)}\ra \nu^{(2)} \ra \nu^{(1)} \ra
\cdots$
 \carre


\vspace{1cm}

\noindent {\bf A3.} The proof of Theorem \ref{T:transform} is more
technical than that of Theorem \ref{T:transform-Rd} and relies on
a series of lemmas stated below.

\begin{lemma}
\label{L:D*} Let $\nu$ be any probability distribution on $[m,M]$,
$0<m \leq M<\infty$ with moments $\mu_i=\int \ml^i \nu(d\ml)$,
$i\in\mathbb{Z}$ ($\mu_0=1$). Then,
\begin{eqnarray}
  \mu_2-\mu_1^2 &\leq& D^* = (M-m)^2/4  \label{D*} \\
    \mu_1\mu_{-1} &\leq& L^* = (M+m)^2/(4mM) \label{L*} \,.
\end{eqnarray}
\end{lemma}

\noindent{\em Proof.} The proof relies on standard results in
experimental design theory, see, e.g., \cite{Fedorov72,Silvey80}.
Consider the two linear regression models
$\eta_1(\mt,\ml)=\mt_0+\mt_1\ml$ and
$\eta_2(\mt,\ml)=\mt_0/\sqrt{\ml}+\mt_1\sqrt{\ml}$, with
$\mt_0,\mt_1$ the model parameters and $\ml$ the design variable,
$\ml\in [m,M]$. $D$-optimum design (approximate theory) aims at
determining a probability measure on $[m,M]$ that maximizes the
determinant of the information matrix associated with a particular
model, here respectively
$$
\Ib_1(\nu) = \left(\begin{array}{cc} \mu_0 & \mu_1 \\ \mu_1 &
\mu_2
\end{array} \right) \ \mbox{ and } \
\Ib_2(\nu) = \left(\begin{array}{cc} \mu_{-1} & \mu_0 \\ \mu_0 &
\mu_1
\end{array} \right) \,.
$$
The function $\log\det \Ib(\nu)$ is concave on the set of
probability measures on $[m,M]$, and its maximum is unique. The
Kiefer-Wolfowitz General Equivalence Theorem \cite{KieferW60}
gives a characterization of the measure $\nu^*$ that maximizes
$\det \Ib_1(\nu)=\mu_2-\mu_1^2$ and $\det
\Ib_2(\nu)=\mu_1\mu_{-1}-1$. In this case it corresponds to the
two point measure, supported at $m$ and $M$, with both masses
equal to $1/2$. Direct calculation gives (\ref{D*},\ref{L*}). One
may notice that (\ref{L*}) corresponds to the Kantorovich
inequality, see \cite{KantorovichA82} and \cite{Luenberger73}, p.\
151. (A full development of this connection is presented in
\cite{PWZa05_LAA}.)

\carre


\begin{lemma}
\label{L:1} Let $\nu$ be any probability distribution on $[m,M]$,
$0<m \leq M<\infty$. Assume that there exists an interval
$\SI\subseteq[m,M]$, $|\SI|\leq \alpha$ and $\nu(\SI)\geq 1-\me$,
$\me\in[0,1]$. Then, $\Var(\nu) \leq \alpha^2/4+2\me M^2$.
\end{lemma}

\noindent{\em Proof.} Define $\mu_1=\int_{[m,M]} \ml \,
\nu(d\ml)$, $\mu_\SI=\int_\SI \ml \, \nu(d\ml)$. Then
$\mu_1=\mu_\SI+\int_{[m,M]\setminus\SI} \ml \, \nu(d\ml)$.
Therefore, $\mu_\SI \leq \mu_1 \leq \mu_\SI + \me M$. We get
\begin{eqnarray*}
\Var(\nu) = \int (\ml-\mu_1)^2 \, \nu(d\ml) &\leq& \int_\SI (\ml-\mu_1)^2 \, \nu(d\ml) + (M-m)^2 \me \\
&& \hspace{-1cm} =\int_\SI (\ml-\mu_\SI)^2 \, \nu(d\ml) +
(\mu_1-\mu_\SI)^2 \nu(\SI) + (M-m)^2 \me \,.
\end{eqnarray*}
Lemma \ref{L:D*} implies $\int_\SI (\ml-\mu_\SI)^2 \,
\nu(d\ml)\leq \alpha^2/4$ and $(\mu_1-\mu_\SI)^2\leq \me^2 M^2$
 gives
$$
\Var(\nu) \leq \alpha^2/4 + \me^2M^2+ M^2\me \leq \alpha^2/4+2\me
M^2 \,.
$$
\carre


\begin{lemma}
\label{L:2} Let $\nu$ be any probability distribution on $[m,M]$,
$0<m \leq M<\infty$. Assume that $\Var(\nu)\leq\me$. Then, there
exist an interval $\SI$ such that $|\SI| \leq \me^{1/4}$ and
$\nu(\SI) \geq 1-4\sqrt{\me}$
\end{lemma}

\noindent{\em Proof.} Take $\SI=[\mu_1-\me^{1/4}/2,
\mu_1+\me^{1/4}/2]$, $\mu_1=\int \ml \, \nu(d\ml)$, and apply the
Chebyshev inequality. \carre


\begin{lemma}
\label{L:3} Let $\nu$ be any distribution on $[m,M]$, $0<m \leq
M<\infty$. Define $\mu_i=\int \ml^i \, \nu(d\ml)$ and
$$
\Mb=\left(
\begin{array}{lll}
\mu_{-1}  & \mu_0 &\mu_1 \\
\mu_0  & \mu_1 &\mu_2 \\
\mu_1  & \mu_2 &\mu_3
\end{array}
\right) \,.
$$
Assume that $L=\mu_{-1}\mu_1>1$ (which, by Jensen's inequality,
holds when $\nu$ is not degenerate at a single point) and $\det
\Mb < \me$. Then, there exist two intervals $\SI_1$ and $\SI_2$
such that
\begin{eqnarray}
&\mbox{(i)}&  |\SI_i| \leq
\frac{(M-m)\me^{1/4}}{m^{9/4}(L-1)^{3/2}} \,, \ i=1,2\,, \
\nu(\SI_1)+\nu(\SI_2) \geq 1-4\sqrt{\me}M^{3/2}\,, \nonumber \\
&\mbox{(ii)}& \max_{x\in\SI_i} |x-\mu_{-1}| > \frac{3(L-1)m^2}{4(M-m)} \,, \ i=1,2 \,, \label{lemma3ii} \\
&\mbox{(iii)}& \mbox{ for }
\me<\me_*=\frac{4(L-1)^8M^8m^{16}}{[32(L-1)^3+M^4(M-m)^2]^2} \,,
\nonumber \\
&& \nu(\SI_i) \geq \frac{m^2(L-1)}{4M^2} \,, \ i=1,2 \,, \label{me*} \\
&& \mbox{ and } \max_{x\in\SI_1,y\in\SI_2} |x-y| > m\sqrt{2(L-1)}
\,. \nonumber
\end{eqnarray}

\end{lemma}

\noindent{\em Proof.} \\
\noindent{\bf (i)} Consider the measure $\nu'$ defined by
$\nu'(\SA)=(1/\mu_{-1}) \int_\SA (1/\ml)\nu(d\ml)$ for any Borel
set $\SA\subset[m,M]$, and denote its moments by
$\mu'_i=(1/\mu_{-1})\int \ml^{i-1} \nu(d\ml)=\mu_{i-1}/\mu_{-1}$.
Note that for any Borel set $\SA$
$$
\frac{1}{M\mu_{-1}} \, \nu(\SA) \leq \nu'(\SA) \leq
\frac{1}{m\mu_{-1}} \, \nu(\SA) \,.
$$
We have
$$
\Mb'=\left(
\begin{array}{lll}
\mu'_0  & \mu'_1 &\mu'_2 \\
\mu'_1  & \mu'_2 &\mu'_3 \\
\mu'_2  & \mu'_3 &\mu'_4
\end{array}
\right) = \Mb/\mu_{-1}
$$
and thus $\det \Mb'=\det \Mb/\mu_{-1}^3$. Also define
$D'=\mu'_2-(\mu'_1)^2$, $a=\sqrt{D'}$,
$b=(\mu'_1\mu'_2-\mu'_3)/\sqrt{D'}$,
$c=a\mu'_2+b\mu'_1=[(\mu'_2)^2-\mu'_1\mu'_3]/\sqrt{D'}$ (note that
$a>0$, $b<0$ and $c<0$) and $\eta=F(\zeta)=a\zeta^2+b\zeta-c$,
with $\zeta$ having the distribution $\nu'$. Direct calculation
gives $E'\{\eta\}=\int \eta(\zeta)\nu'(d\zeta)=0$ and
$\Var'(\eta)=E'\{\eta^2\}-(E'\{\eta\})^2=\det \Mb'$, so that $\det
\Mb < \me$ implies $\Var'(\eta) < \me'=\me/\mu_{-1}^3$. From Lemma
\ref{L:2}, the interval $\SI=[-(\me')^{1/4}/2, (\me')^{1/4}/2]$ is
such that $\Pr\{\eta\in \SI\} \geq 1-4\sqrt{\me'}$. Also, from the
mean-value theorem, there exist $\ml_1<\ml_2$ such that
$\ml_i\in[m,M]$ and $a\ml_i^2+b\ml_i-c=0$, $i=1,2$. Direct
calculation gives
$F(\mu'_1)=F(\mu_{-1})=a(\mu'_1)^2+b\mu'_1-c=-(D')^{3/2}$, and
thus
$$
m \leq \ml_1 < 1/\mu_{-1} < \ml_2 \leq M \,.
$$
Take $\beta=(M-m)(\me')^{1/4}/[2(D')^{3/2}]$, we get
\begin{eqnarray*}
  F(\ml_1+\beta) &<& -(\me')^{1/4}/2\,, F(\ml_1-\beta)>(\me')^{1/4}/2\,, \\
  F(\ml_2+\beta) &>& (\me')^{1/4}/2\,, F(\ml_2-\beta)<-(\me')^{1/4}/2\,,
\end{eqnarray*}
and $\nu(\SI_1)+\nu(\SI_2) \geq 1-4\sqrt{\me}M^{3/2}$ when
$\SI_i=[\ml_i-\beta,\ml_i+\beta]$, $i=1,2$, with
$|\SI_i|=2\beta=(M-m)\me^{1/4}\mu_{-1}^{9/4}/(L-1)^{3/2}\leq
(M-m)\me^{1/4}/[m^{9/4}/(L-1)^{3/2}]$.

\vspace{0.3cm} \noindent{\bf (ii)} Define $y_1=\mu'_1-\ml_1$,
$y_2=\ml_2-\mu'_1$, so that $\max_{x\in\SI_1} |x-\mu_{-1}|
> y_1$ and $\max_{x\in\SI_2} |x-\mu_{-1}|
> y_2$. We have $F(\ml)=a(\ml-\ml_1)(\ml-\ml_2)$ and thus
$y_1y_2=-F(\mu'_1)/a=D'$. Also, $|y_2-y_1|<y_1+y_2\leq M-m$, so
that $D' > y_i(y_i+M-m)$, $i=1,2$, and thus
$$
y_i > \frac{M-m}{2} \, \left[\sqrt{1+\frac{4D'}{(M-m)^2}} -1
\right]
> \frac{D'}{M-m} \, \left(1-\frac{D'}{(M-m)^2}\right) \,, \ i=1,2
\,.
$$
Lemma \ref{L:D*} gives $D'<(M-m)^2/4$, so that
$$
y_i>3D'/[4(M-m)]>3(L-1)m^2/[4(M-m)]\,, \ i=1,2 \,.
$$

\noindent{\bf (iii)} Define $\gamma=\nu'(\SI_2)$, part (i) implies
$\nu'(\SI_1) > 1-4\sqrt{\me'}-\gamma$, and from Lemma~\ref{L:1}
$$
D'\leq \frac{(M-m)^2\sqrt{\me'}}{4(D')^3} + 2
(4\sqrt{\me'}+\gamma)M^2 \,,
$$
which gives
$$
\gamma \geq \frac{D'}{2M^2} - \sqrt{\me'} \, \left[
\frac{(M-m)^2}{8(D')^3M^2}+4 \right] \,,
$$
and thus $\gamma \geq D'/(4M^2)>(L-1)m^2/[4M^2]$ for
$\me<\me_*<[4(D')^8]/[(M-m)^2+32(D')^3M^2]^2$, see (\ref{me*}).

Define now $\Delta=\max_{x\in\SI_1,y\in\SI_2} |x-y|$. Lemma
\ref{L:1} gives $D'\leq \Delta^2/4+8\sqrt{\me'}M^2$, which implies
$\Delta^2 \geq 4D'-32M^2\sqrt{\me'}$. Since $\me<\me_*$ implies
$\sqrt{\me'}<D'/(16M^2)$, we get $\Delta^2>2(L-1)m^2$. \carre


\vspace{0.5cm} \noindent{\em Proof of Theorem \ref{T:transform}.}
The proof follows the same lines as that of Theorem
\ref{T:transform-Rd} and is divided into four parts. In (i), we
construct sequences of intervals $\SL_k=[m_k,m_k+\delta]$ and
$\SR_k=[M_k-\delta,M_k]$ in which the measure $\nu_k$ will tend to
concentrate. In (ii) we prove that $\SR_k \cap \SR_{k+1} \neq
\emptyset$ and in (iii) that the sequence $M_k$ is non-decreasing.
Finally, the limiting behaviour of $\nu_k$ is derived in (iv).

\vspace{0.3cm} \noindent{\bf (i)} We have seen in
Section~\ref{S:monotonicity} that $\det \Mb_k\ra 0$ as
$k\ra\infty$, with $\Mb_k$ given by (\ref{Mk}). Therefore, given $
\me$, $\exists K_\me$ such that $\forall k>K_\me$, $\det
\Mb_k<\me$. Define $L_k=\mu_1^k\mu_{-1}^k$ and note that $L_k>1$
because no $\nu_k$ is degenerate at a single point. Using
Lemma~\ref{L:3}, for $\me$ small enough, for any $k>K_\me$ there
exist two intervals $\SI_1^k$, $\SI_2^k$, with width at most
$$
\delta = \delta(\me)= \frac{ (M-m)\me^{1/4}}{m^{9/4}(L_0-1)^{3/2}}
\,,
$$
and such that $\nu_k(\SI_1^k)+\nu_k(\SI_2^k)\geq
1-4\sqrt{\me}M^{3/2}$, $\nu_k(\SI_1^k)\geq m^2(L_k-1)/(4M^2)$,
$\nu_k(\SI_2^k)\geq m^2(L_k-1)/(4M^2)$. Also,
$\max_{x\in\SI_1^k,y\in\SI_2^k}|x-y| \geq m\sqrt{2(L_k-1)}$.
Without any loss of generality, assume that $\SI_1^k$ is the
interval on the left. Define $\SL(x)=[x,x+\delta]$,
$\SR(x)=[x-\delta,x]$,
\begin{eqnarray*}
\SX_L^k &=& \Arg \max_x \{\nu_k[\SL(x)]\,, \, \SL(x)\cap\SI_1^k \neq \emptyset \} \,, \\
\SX_R^k &=& \Arg \max_x \{\nu_k[\SR(x)]\,, \, \SR(x)\cap\SI_2^k
\neq \emptyset \} \,,
\end{eqnarray*}
and $m_k=\min \SX_L^k$, $M_k=\max \SX_R^k$, $\SL_k=\SL(m_k)$,
$\SR_k=\SR(M_k)$; that is, $M_k$ is the right endpoint of an
interval $\SR_k$, intersecting $\SI_2^k$, with maximum measure,
and similarly for $m_k$ and $\SL_k$. Note that
$\nu_k(\SL_k)+\nu_k(\SR_k)\geq 1-4\sqrt{\me}M^{3/2}$,
$\nu_k(\SL_k)\geq m^2(L_k-1)/(4M^2)$ and $\nu_k(\SR_k)\geq
m^2(L_k-1)/(4M^2)$. The situation is the same for the two
sequences of intervals $(\SL_k)$ and $(\SR_k)$, and we concentrate
on $(\SR_k)$ in the rest of the proof.

\vspace{0.3cm}
\noindent{\bf (ii)} We show now that $\SR_k \cap
\SR_{k+1} \neq \emptyset$. Again for $\me$ small enough $\mu_1^k
\notin \SR_k$ and $\ml-\mu_1^k \geq M_k-\delta-\mu_1^k$ on $\SR_k$
so that
\begin{eqnarray*}
\nu_{k+1}(\SR_k) &=& \int_{\SR_k} \frac{(\ml-\mu_1^k)^2}{D_k} \,
\nu_k(d\ml) \geq \frac{\nu_k(\SR_k)}{D_k}
(M_k-\delta-\mu_1^k)^2 \\
&& \geq \frac{m^2(L_k-1)}{4M^2D^*} (M_k-\delta-\mu_1^k)^2
\end{eqnarray*}
with $D^*$ the maximum possible value of $D_k$, $D^*=(M-m)^2/4$,
see Lemma \ref{L:D*}. By construction, $\max_{x\in\SI_k^2}
|x-\mu_1^k| \leq M_k+\delta-\mu_1^k$, and thus, from Lemma
\ref{L:3},
\begin{equation}
M_k-\mu_1^k+\delta > \frac{3m^2(L_k-1)}{4(M-m)} \geq
\frac{3m^2(L_0-1)}{4(M-m)} = C \,. \label{C**}
\end{equation}
Choosing $\me$ such that $\delta<C/4$ gives
$M_k-\delta-\mu_1^k>C/2$ and thus
$$
\nu_{k+1}(\SR_k) > \frac{m^2(L_k-1)}{4M^2D^*} \, \frac{C^2}{4}
\geq \nu_R^*=\frac{9m^6(L_0-1)^3}{16M^2(M-m)^4} \,.
$$
Choosing now $\me$ such that $4\sqrt{\me}M^{3/2} < \nu_R^*$ we
obtain $\SR_k \cap \SR_{k+1} \neq \emptyset$ for any $k>K_\me$.

\vspace{0.3cm}
\noindent{\bf (iii)} We prove now that the sequence
$(M_k)$ is not decreasing starting at some $K_\me$ for $\me$ small
enough. Take $k>K_\me$ and assume that $M_{k+1}=M_k-\beta$,
$\beta>0$. Then note that  $\beta<\delta$ since $\SR_k \cap
\SR_{k+1} \neq \emptyset$ by (ii) above. Consider the difference
$\nu_{k+1}(\SR_k)-\nu_{k+1}(\SR_{k+1})=\nu_{k+1}([M_k-\beta,M_k])-\nu_{k+1}([M_k-\delta-\beta,M_k-\delta])$.
Assume first that $\nu_{k+1}([M_k-\delta-\beta,M_k-\delta])=0$,
then $\nu_{k+1}(\SR_k)>\nu_{k+1}(\SR_{k+1})$, which is impossible
by construction. We can thus consider the following ratio
\begin{eqnarray*}
\frac{\nu_{k+1}([M_k-\beta,M_k])}{\nu_{k+1}([M_k-\delta-\beta,M_k-\delta])}
&=&
\frac{\int_{M_k-\beta}^{M_k} (\ml-\mu_1^k)^2 \, \nu_k(d\ml) }{\int_{M_k-\delta-\beta}^{M_k-\delta} (\ml-\mu_1^k)^2 \, \nu_k(d\ml)} \\
&& \geq \frac{(M_k-\beta-\mu_1^k)^2}{(M_k-\delta-\mu_1^k)^2} \,
\frac{\nu_{k}([M_k-\beta,M_k])}{\nu_{k}([M_k-\delta-\beta,M_k-\delta])}
\,.
\end{eqnarray*}
Since $M_k-\delta-\mu_1^k \geq C-2\delta \geq 2\delta$ for
$C>4\delta$, see (\ref{C**}), and $\beta<\delta$,
$(M_k-\beta-\mu_1^k)^2 > (M_k-\delta-\mu_1^k)^2$. Also, by
construction,
\begin{eqnarray*}
  0 &\leq& \nu_{k}([M_k-\beta,M_k])
-\nu_{k}([M_k-\delta-\beta,M_k-\delta]) \\
  && = \nu_{k}([M_k-\beta,M_k])-\nu_{k}([M_k-\delta-\beta,M_k-\delta])\,.
\end{eqnarray*}
This gives
$$
\frac{\nu_{k+1}([M_k-\beta,M_k])}{\nu_{k+1}([M_k-\delta-\beta,M_k-\delta])}
> 1 \,.
$$
Therefore, $\beta>0$ leads to
$\nu_{k+1}(\SR_k)>\nu_{k+1}(\SR_{k+1})$, which is impossible. We
thus obtain $M_{k+1}\geq M_k$ for $k>K_\me$.

\vspace{0.3cm} \noindent{\bf (iv)} Since the sequence $(M_k)$ is
non-decreasing and bounded from above (by $M$), it has a limit
$M_*\geq M$. The same is true for $m_k$, and $m_k\ra m_*$ as
$k\ra\infty$. We have thus proved that for any $\delta$ small
enough and any $k$ larger than some $K_\delta$,
$$
\nu_k([M_*-\delta,M_*])+\nu_k([m_*,m_*+\delta]) \geq 1 -
\frac{4M^{3/2}m^{9/2} (L_0-1)^3\delta^2 }{(M-m)^2} \,.
$$
Assume that $M_*<M$. This would imply $\nu_k([M-\delta,M]) \ra 0$
as $k\ra\infty$ for $\delta<M-M_*$. On the other hand,
$$
\frac{\nu_{k+1}([M-\delta,M])}{\nu_{k+1}([M_*-\delta,M_*])} >
\frac{\nu_{k}([M-\delta,M])}{\nu_{k}([M_*-\delta,M_*])} \,,
$$
which leads to a contradiction since
$\nu_{k}([M-\delta,M])/\nu_{k}([M_*-\delta,M_*])$ is then
increasing and $\nu_{k}([M_*-\delta,M_*])$ is bounded from below.
Therefore, $M_*=M$, and similarly $m_*=m$, with, for $\delta$
small enough and any $k$ larger than some $K_\delta$,
$\nu_k([m+\delta,M-\delta]) < 4M^{3/2}m^{9/2} (L_0-1)^3\delta^2
/(M-m)^2$. Finally, from Helly's Theorem, see \cite{Shiryaev96},
p.\ 319, from the sequence $(\nu_k)$ we can extract a  subsequence
$(\nu_{k_i})$ that is weakly convergent, and from the result above
the associated limit has necessarily the form $\nu_p^*$, where
$\nu_p^*$ is the discrete measure concentrated on the two points
$m$, $M$, with $\nu_p^*(m)=p$, $\nu_p^*(M)=1-p$. Since $L_{k_i}$
converges to some $L$, $\nu_p^*$ is such that the associated value
of $\mu_1\mu_{-1}$ is equal to $L$, which only leaves two
possibilities for $p$ (and $1-p$):
$$
p= \frac{1}{2} \pm \frac{\rho+1}{\rho-1} \sqrt{\frac{1}{4}-
\frac{\rho L}{(\rho+1)^2}}
$$
where $\rho=M/m$. Applying the transformation $T$, we get
$\nu_{k_i+1}=T(\nu_{k_i}) \ra T(\nu_p^*)=\nu_{1-p}^*$. \carre

\vspace{1cm}

\noindent  {\bf A4.} {\em Proof of Theorem \ref{T:instability}.}

\noindent{\bf (i)} It is straightforward to check that
$T^2(\nu_p^*)=\nu_p^*$, $\forall p\in(0,1)$.

\vspace{0.3cm}
\noindent{\bf (ii)} We assume that $\SS_A$ is not
reduced to $\{m,M\}$ (otherwise $\SI_u=\emptyset$). We have
$\nu_{k+2}(d\ml) = H(\nu_k,\ml) \nu_k(d\ml)$, with
\begin{equation}\label{Hnuk}
H(\nu_k,\ml)=\frac{(\ml-\mu_1^k)^2(\ml-\mu_1^{k+1})^2}{D_kD_{k+1}}
\end{equation}
see (\ref{nuk+1}), with $\mu_1^k$, $D_k$ defined as in Theorem
\ref{T:transform}. For $\nu_k=\nu_p^*$, it gives
\begin{equation}\label{Hnup*}
    H(\nu_p^*,\ml)=\frac{[M(1-p)+m p-\ml]^2[M
p+m(1-p)-\ml]^2}{p^2(1-p)^2(M-m)^4} \,.
\end{equation}
One can then check that for any $p\in\SI_u$, $\max_{\ml\in\SS_A}
H(\nu_p^*,\ml)=H(\nu_p^*,\ml^*)>1$, with $\ml^*=\min_{\ml\in\SS_A}
s(\ml)$. Therefore, for any $p\in\SI_u$, one can choose $\me$
small enough, such that $d(\nu_k,\nu_p^*) < \me$ implies
$\nu_{k+2}([a,b]) > K_p \nu_k([a,b])$, for some $K_p>1$ and some
$a,b$ such that $m+\me<a<b<M-\me$ and
$[a,b]\cap\SS_A\neq\emptyset$. For any $\alpha>0$, $\alpha<1-p$,
take an initial measure $\nu_0$ putting mass $p$ at $m$,
$1-p-\alpha$ at $M$ and $\alpha$ in the interval $[a,b]$. It
satisfies $d(\nu_0,\nu_p^*)<\alpha$, and, for any $m$, either
$d(\nu_{2m},\nu_p^*)>\me$ or $\nu_{2(m+1)}([a,b])> K_p
\nu_{2m}([a,b])$. The later case gives $\nu_{2m}([a,b])> 2\me$,
and thus $d(\nu_{2m},\nu_p^*)>\me$, as soon as
$m>\log(2\me/\alpha)/\log(K_p)$, which shows that $\nu_p^*$ is
unstable.

\vspace{0.3cm}
\noindent{\bf (iii)} Part (a) concerns the case
where a spectral gap is present, with point mass at $m$ and $M$.
The proof for the general situation is more technical and is
sketched in part (b).

\vspace{0.3cm} {\bf (a)} Assume that the measure $\nu_0$ has a
spectral gap: $\nu_0=0$ on $(m,m+s)$ and $(M-s,M)$ for some $s>0$.
Take $\gamma<s$ and assume that $d(\nu_0,\nu_p^*)<\alpha<\gamma$
with $p\in\SI_s$. The arguments go as follows. First we bound
$\nu_{2}\{(m+\gamma,M-\gamma]\}$ by $2K_0\alpha$ for some $K_0<1$,
then we bound $\nu_{2}\{(M-\gamma,M]\}$ by $1-p+K_1\alpha$ for
some $K_1<\infty$. We show that $d(\nu_2,\nu_{p_2}^*)<K_0\alpha$
for some $p_2$ such that $|p_2-p|<(K_0+K_1)\alpha$. Stability will
then follow by an induction argument.

The maximum value of $H(\nu_0,\ml)$ for $\ml$ varying in
$[m+\gamma,M-\gamma]$ may be reached for some
$\ml^*\in(\mu_1^0,\mu_1^{1})$ or at one of the two points
$m+\gamma$, $M-\gamma$. Now, for $\alpha$ small enough
$H(\nu_0,\ml)$ will be close to $H(\nu_p^*,\ml)$ given by
(\ref{Hnup*}), and $p\in\SI_s$ implies
\begin{equation}\label{*2}
    \max_{\ml\in\SS_A\cap(\mu_1^0,\mu_1^{1})} H(\nu_0,\ml)<1 \,.
\end{equation}
Consider the function $H(\nu_0,\ml)$ at $\ml=M-\gamma$. We can
write
\begin{equation}\label{**}
    H(\nu_0,M-\gamma)=H(\nu_p^*,M)-\gamma \,
    \frac{dH(\nu_p^*,\ml)}{d\ml}_{|\ml=M} + F_H(\nu_p^*;\nu_0,M) +
    \SO(\gamma^2) \,,
\end{equation}
with $F_H(\nu_p^*;\nu_0,M)$ the directional derivative of
$H(\nu,M)$ at $\nu_p^*$ in the direction $\nu_0$,
$$
F_H(\nu_p^*;\nu_0,M) = \lim_{\beta\ra 0^+} \frac{
H[(1-\beta)\nu_p^*+\beta\nu_0,M]-H(\nu_p^*,M)}{\beta} \,.
$$
Define $F_H(\nu_p^*,x,\ml)=F_H(\nu_p^*;\delta_{x},\ml)$ with
$\delta_x$ the delta measure supported at $x$. We have
$$
F_H(\nu_p^*;\nu_0,M) = \int_m^M F_H(\nu_p^*,x,M) \nu_k(dx) \,,
$$
which we decompose in three parts:
\begin{eqnarray*}
  F_H(\nu_p^*;\nu_0,M) &=& \int_m^{m+\gamma} F_H(\nu_p^*,x,M)
\nu_0(dx) + \int_{m+\gamma}^{M-\gamma} F_H(\nu_p^*,x,M) \nu_0(dx) \\
  && + \int_{M-\gamma}^M F_H(\nu_p^*,x,M) \nu_0(dx) \,.
\end{eqnarray*}
Direct calculation gives
$$
F_H(\nu_p^*,x,M) =
\frac{(x-m)^2(M-x)[x-m+(2p-1)(M-m)]}{p^2(1-p)^2(M-m)^4}
$$
so that $F_H(\nu_p^*,m,M)=F_H(\nu_p^*,M,M)=0$ and
$F_H(\nu_p^*;\nu_0,M) < F^* \nu_0\{(m+\gamma,M-\gamma]\}$ with
$F^*=\max_{p\in\SI_s,\,x\in[m,M]} F_H(\nu_p^*,x,M)<\infty$. Also,
$d(\nu_0,\nu_p^*)<\alpha$ implies
$\nu_0\{(m+\gamma,M-\gamma]\}=\nu_0\{(m+\alpha,M-\alpha]\}<2\alpha$,
so that $F_H(\nu_p^*;\nu_0,M) < 2\alpha F^*$. Now,
$$
H(\nu_p^*,M)=1\,, \ \frac{d
H(\nu_p^*,\ml)}{d\ml}_{|\ml=M}=\frac{2}{p(1-p)(M-m)} \,,
$$
which, together with (\ref{**}) gives for $\gamma$ small enough
$$
  H(\nu_0,M-\gamma) < 1+2\alpha F^* - \frac{\gamma}{p(1-p)(M-m)}
$$
and thus
$$
  H(\nu_0,M-\gamma) < 1 - \frac{\gamma}{2p(1-p)(M-m)}
$$
for $\alpha<\gamma/[4F^*p(1-p)(M-m)]$.

The situation is similar at $m+\gamma$. Together with (\ref{*2})
this implies for $\alpha$ small enough
$$
\max_{\ml\in\SS_A\cap[m+\gamma,M-\gamma]} H(\nu_0,\ml) < K_0 < 1
$$
and therefore,
\begin{equation}\label{*3}
    \nu_{2}\{(m+\gamma,M-\gamma]\} < 2K_0\alpha
\end{equation}
with $K_0<1$ not depending on $\alpha$.

Consider now the interval $(M-\gamma,M]$. We have
$$
\nu_{2}\{(M-\gamma,M]\} = \nu_{2}(M) = H(\nu_0,M)
\nu_{0}\{(M-\gamma,M]\} = H(\nu_0,M) \nu_0(M) \,,
$$
with $d(\nu_0,\nu_p^*)<\alpha$ implying $\nu_{0}(M)<1-p+\alpha$,
and
$$
H(\nu_0,M)=H(\nu_p^*,M)+F_H(\nu_p^*;\nu_0,M)+\SO(\alpha^2) <
1+2\alpha F^* + \SO(\alpha^2) \,.
$$
This gives for $\alpha$ small enough
$$
\nu_{2}\{(M-\gamma,M]\} < 1-p+K_1 \alpha
$$
for some $K_1<\infty$. Similarly, $\nu_{2}\{[m,m+\gamma]\} < p+K_1
\alpha$.

Define $p_0=p$, $p_{2}=[\nu_{2}(m)-\nu_{2}(M)+1]/2$,
$\alpha_0=\alpha$. We obtain
$$
-K_0\alpha_0 < \nu_{2}(m)-p_{2} < 0\,, \ -K_0\alpha_0 <
\nu_{2}(M)-(1-p_{2}) < 0
$$
which together with (\ref{*3}) implies
$$
d(\nu_{2},\nu_p^*) < \alpha_{2}=K_0 \alpha_0 \,.
$$
Moreover, $|p_{2}-p_0|<(K_0+K_1)\alpha_0$.

For $\alpha$ small enough, $p_2\in\SI_s$ and we can then repeat
the same arguments. This gives for any $m$
$$
d(\nu_{2m},\nu_{p_{2m}}^*) < \alpha_{2m}=K_0^m \alpha
$$
with
$$
|\, p_{2m}-p\,| < (K_0+K_1) \sum_{i=0}^{m-1} \alpha_{2i} =
(K_0+K_1) \frac{1-K_0^m}{1-K_0} \, \alpha < \frac{K_0+K_1}{1-K_0}
\, \alpha
$$
and $p_{2m}\in\SI_s$, for $\alpha$ small enough. For any
$p\in\SI_s$ and any $\me>0$, taking $\nu_0$ such that
$d(\nu_0,\nu_p^*)<\alpha$ with $\alpha$ small enough thus implies
$d(\nu_{2m},\nu_p^*)<\me$ for any $m$, and $\nu_p^*$ is thus
stable.

\vspace{0.3cm}
{\bf (b)} Consider now the general situation. The
proof follows the same lines as in case (a), but more
technicalities are required since we need to consider measures of
intervals of the form $[m,m+\gamma]$ and $(M-\gamma,M]$, with
$\gamma$ decreasing in a suitable way as the number of iterations
of the mapping $T^2$ increases.

Assume that
\begin{eqnarray*}
  && \nu_{2k}\{(m+\gamma_{2k},M-\gamma_{2k}]\} < 2 \alpha_{2k} \,, \\
  && \nu_{2k}\{[m,m+\gamma_{2k}]\} < p_{2k}+ \alpha_{2k}\,, \\
  && \nu_{2k}\{(M-\gamma_{2k},M]\} < 1-p_{2k} + \alpha_{2k} \,.
\end{eqnarray*}
for some $p_{2k}\in\SI_s$ and some $\alpha_{2k}$, $\gamma_{2k}$.
Note that it implies $d(\nu_{2k},\nu_{p_{2k}}^*)<\gamma_{2k}$ and
that for $k=0$, $\alpha_{0}$, $\gamma_{0}$ can be chosen
arbitrarily small, with $d(\nu_0,\nu_p^*)<\alpha_0$ for some
$p\in\SI_s$.

Consider one application of the mapping $T^2$ at a generic
iteration $k$. We can write
$H(\nu_{2k},M)=H(\nu_{p_{2k}}^*,M)+F_H(\nu_{p_{2k}}^*;\nu_{2k},M)+\SO(\gamma_{2k}^2)$
with
\begin{eqnarray*}
  F_H(\nu_{p_{2k}}^*;\nu_{2k},M) &=& \int_m^{m+\gamma_{2k}}
F_H(\nu_{p_{2k}}^*,x,M) \nu_{2k}(dx) \\
&& + \int_{m+\gamma_{2k}}^{M-\gamma_{2k}} F_H(\nu_{p_{2k}}^*,x,M)
\nu_{2k}(dx)  \\
  && + \int_{M-\gamma_{2k}}^M F_H(\nu_{p_{2k}}^*,x,M) \nu_{2k}(dx) \,.
\end{eqnarray*}
The first integral term is of the order $\SO(\gamma_{2k}^2)$
(since $F_H(\nu_{p_{2k}}^*,m,M)=0$ and
$dF_H(\nu_{p_{2k}}^*,z,M)/dz_{|z=m}=0$), the second is bounded by
$2\alpha_{2k}F^*+\SO(\gamma_{2k}^2)$, as in case (a). For the
third term, for which $x$ is close to $M$, we can use the linear
approximation
\begin{eqnarray*}
  F_H(\nu_{p_{2k}}^*,x,M) &=& (x-M) \frac{d
F_H(\nu_{p_{2k}}^*,z,M)}{dz}_{|z=M} + \SO(\gamma_{2k}^2) \\
  &=& \frac{-2(x-M)}{p_{2k}(1-p_{2k})^2(M-m)} + \SO(\gamma_{2k}^2)
\end{eqnarray*}
which gives
$$
\int_{M-\gamma_{2k}}^M F_H(\nu_{p_{2k}}^*,x,M) \nu_{2k}(dx) =
\frac{2}{p_{2k}(1-p_{2k})^2(M-m)} I_{2k}(M) + \SO(\gamma_{2k}^2)
$$
where $I_{2k}(M)=\int_0^{\gamma_{2k}} z \nu_{2k}'(dz)$ with
$\nu_{2k}'$ the measure obtained after applying the transformation
$x \mapsto z=M-x$. We have thus obtained
\begin{equation}\label{boundH}
H(\nu_{2k},M) < 1 + 2\alpha_{2k}F^* +
\frac{2I_{2k}(M)}{p_{2k}(1-p_{2k})^2(M-m)} + \SO(\gamma_{2k}^2)
\,.
\end{equation}
Consider now the behavior of $I_{2k}(M)$ as $k$ increases. We
assume that $\nu_{2k}$ remains in some neighborhood $\SV(p)$ of
$\nu_p^*$, which we shall be able to guarantee afterwards. Define
$A_{2k}(M)=I_{2k}(M)[\int_0^{\gamma_{2k}} \nu_{2k}'(dz)]^{-1}$. It
satisfies $I_{2k}(M)<A_{2k}(M)<\gamma_{2k}$. Also,
$\gamma_{2(k+1)}<\gamma_{2k}$ implies
$$
A_{2(k+1)}(M)=\frac{\int_0^{\gamma_{2(k+1)}} z H(\nu_{2k},M-z)
\nu_{2k}'(dz)} {\int_0^{\gamma_{2(k+1)}} H(\nu_{2k},M-z)
\nu_{2k}'(dz)} < \frac{\int_0^{\gamma_{2k}} z H(\nu_{2k},M-z)
\nu_{2k}'(dz)} {\int_0^{\gamma_{2k}} H(\nu_{2k},M-z)
\nu_{2k}'(dz)} \,,
$$
and, since $H(\nu_{2k},M-z)$ decreases for $z$ close to zero,
$$
A_{2(k+1)}(M) < \frac{\int_0^{\gamma_{2k}} z
\frac{H(\nu_{2k},M-z)}{H(\nu_{2k},M)} \nu_{2k}'(dz)}
{\int_0^{\gamma_{2k}} \nu_{2k}'(dz)} \,.
$$
We can bound the speed of decrease of $H(\nu_{2k},M-z)$:
$H(\nu,M-z)/H(\nu,M)<1-az$ for some $a>0$, any $z$ in
$[0,\gamma_0]$ and any $\nu\in\SV(p)$. This gives
$$
A_{2(k+1)}(M) < \frac{\int_0^{\gamma_{2k}} z (1-az) \nu_{2k}'(dz)}
{\int_0^{\gamma_{2k}} \nu_{2k}'(dz)} \,.
$$
Repeating the same arguments we get for any $n>0$,
$$
A_{2(k+n)}(M) < \bar{A}_{2(k+n)}(M) = \frac{\int_0^{\gamma_{2k}} z
(1-az)^n \nu_{2k}'(dz)} {\int_0^{\gamma_{2k}} \nu_{2k}'(dz)} \,,
$$
with $\bar{A}_{2(k+n)}(M)$ decreasing with $n$. Direct calculation
gives $\sum_{n=0}^\infty \bar{A}_{2(k+n)}(M) =1/a$, and therefore
$I_{2k}(M)<\bar{A}_{2k}(M)=\So(1/k)$.

Similarly to case (a), we can write
$$
H(\nu_{2k},M-\gamma_{2(k+1)})=H(\nu_{2k},M) -
\frac{2\gamma_{2(k+1)}}{p_{2k}(1-p_{2k})(M-m)} +
\SO(\gamma_{2k}^2) \,,
$$
with $H(\nu_{2k},M)$ bounded by (\ref{boundH}). Assume that
$\gamma_{2k}$ is such that $\bar{A}_{2k}(M)=\So(\gamma_{2k})$ and
$\alpha_{2k}=\So(\gamma_{2k})$. We obtain for $p_{2k}$ close
enough to $p$
\begin{equation}\label{beta}
    H(\nu_{2k},M-\gamma_{2(k+1)}) < \beta_{2(k+1)} =
1-\frac{\gamma_{2(k+1)}}{p(1-p)(M-m)} \,.
\end{equation}
We thus get the following bounds on the measure of subintervals of
interest at the next iteration:
\begin{equation}\label{*4}
    \nu_{2(k+1)}\{(m+\gamma_{2(k+1)},M-\gamma_{2(k+1)}]\} < 2 \max\{
\beta_{2(k+1)}, K_0\} \alpha_{2k}
\end{equation}
where $K_0= \max_{\nu_{2k}\in\SV(p)}
\max_{\ml\in\SS_A\cap(\mu_1^{2k},\mu_1^{2k+1})} H(\nu_{2k},\ml)$,
and $K_0<1$ for $p$ in $\SI_s$ and $\SV(p)$ small enough, see part
(a);
\begin{eqnarray*}
  \nu_{2(k+1)}\{(M-\gamma_{2(k+1)},M]\} &<&
\nu_{2(k+1)}\{(M-\gamma_{2k},M]\} \\
&& \hspace{-4cm} < H(\nu_{2k},M) \nu_{2k}
\{(M-\gamma_{2k},M]\} \\
&& \hspace{-4cm} <  \left[ 1 + 2\alpha_{2k}F^* +
\frac{2\bar{A}_{2k}(M)}{p_{2k}(1-p_{2k})^2(M-m)} +
\SO(\gamma_{2k}^2) \right] \nu_{2k} \{(M-\gamma_{2k},M]\} \\
&& < \nu_{2k} \{(M-\gamma_{2k},M]\} + B \alpha_{2k} + C
\bar{A}_{2k}(M) + D \gamma_{2k}^2
\end{eqnarray*}
for some $B,C,D<\infty$. Similarly, we obtain
$$
\nu_{2(k+1)}\{[m,m+\gamma_{2(k+1)}]\} < \nu_{2k}
\{[m,m+\gamma_{2k}]\} + B \alpha_{2k} + C \bar{A}_{2k}(m) + D
\gamma_{2k}^2
$$
where $\bar{A}_{2k}(m)$ is defined similarly to $\bar{A}_{2k}(M)$.
Define $p_{2(k+1)}$ as
$$
p_{2(k+1)} = \frac{ \nu_{2(k+1)}\{[m,m+\gamma_{2(k+1)}]\} -
\nu_{2(k+1)}\{(M-\gamma_{2(k+1)},M]\} +1}{2} \,,
$$
it gives
\begin{eqnarray*}
  && 0 < p_{2(k+1)} - \nu_{2(k+1)}\{[m,m+\gamma_{2(k+1)}]\} <  \max\{
\beta_{2(k+1)}, K_0\} \alpha_{2k} \,, \\
  && 0 < 1-p_{2(k+1)} - \nu_{2(k+1)}\{(M-\gamma_{2(k+1)},M]\} <  \max\{
\beta_{2(k+1)}, K_0\} \alpha_{2k} \,.
\end{eqnarray*}
Together with (\ref{*4}) it implies
$d(\nu_{2(k+1)},\nu_{p_{2(k+1)}}^*)<\gamma_{2(k+1)}< \gamma_{2k}$,
with
$$
|p_{2(k+1)}-p_{2k}| < \Delta_{2k}=[B+1+\max\{ \beta_{2(k+1)},
K_0\}] \alpha_{2k} + C \bar{A}'_{2k} +  D \gamma_{2k}^2 \,,
$$
where $\bar{A}'_{2k}=\max\{\bar{A}_{2k}(m),\bar{A}_{2k}(M)\}$ and
$\sum_k \bar{A}'_{2k}<\infty$.

Define $\alpha_{2(k+1)}=\max\{ \beta_{2(k+1)}, K_0\} \alpha_{2k}$
and take $\gamma_{2k}=1/k^q$ with $q<1$, so that
$\bar{A}'_{2k}=\So(\gamma_{2k})$. From the definition of
$\beta_{2(k+1)}$, see (\ref{beta}), $\sum_{k} \alpha_{2k} <
\infty$ and $\alpha_{2k}=\So(\gamma_{2k})$. Since $\sum_k
\bar{A}'_{2k}<\infty$, taking $q>1/2$ in the definition of
$\gamma_{2k}$ ensures $\sum_{k} \Delta_{2k}<\infty$. We can repeat
the same argument, and
$d(\nu_{2(k+n)},\nu_{p_{2(k+n)}}^*)<\gamma_{2(k+n)}$ which tends
to zero as $n$ increases, with $|p_{2(k+n)}-p_{2k}|$ remaining
finite. $\nu_{2(k+n)}$ thus remains in some neighborhood $\SV(p)$
of $\nu_p^*$ for any $n$, and $\SV(p)$ can be made arbitrarily
small by choosing $\alpha_{0}$ and $\gamma_{0}$ small enough.
\carre

\vspace{1cm}

\noindent {\bf A5.} {\em Proof of Theorem
\ref{T:rates-equivalent}.} Assume that $x_0$ is such that for some
$k\geq 0$, $\|g_{k+1}\|=0$ with $\|g_i\|>0$ for all $i\leq k$
(that is, $x_{k+1}=x^*$ and $x_i\neq x^*$ for $i\leq k$). This
implies $R_k(W)=0$ for any $W$, and therefore
$R(W,x_0,x^*)=R(x_0,x^*)=0$.

Assume now that $\|g_k\|>0$ for all $k$. Consider
$$
V_n=\left[\prod_{k=0}^{n-1} R_k(W) \right]^{1/n} =
\left[\prod_{k=0}^{n-1} \frac{(Wg_{k+1},g_{k+1})}{(Wg_k,g_k)}
\right]^{1/n} = \left[ \frac{(Wg_n,g_n)}{(Wg_0,g_0)}\right]^{1/n}
\,.
$$
We have,
$$
\forall z\in\SH\,, \ c \|z\|^2 \leq (Wz,z) \leq C \|z\|^2 \,,
$$
and thus
$$
(c/C)^{1/n} \left[ \frac{(g_n,g_n)}{(g_0,g_0)}\right]^{1/n} \leq
V_n \leq (C/c)^{1/n} \left[
\frac{(g_n,g_n)}{(g_0,g_0)}\right]^{1/n} \,.
$$
Since $(c/C)^{1/n}\ra 1$ and $(C/c)^{1/n}\ra 1$ as $n\ra\infty$,
$\liminf_{n\ra\infty} V_n$ and $\limsup_{n\ra\infty}$ $V_n$ do not
depend on $W$. Take $W=P(A)$; it gives $R_k(W)=r_k=1-1/L_k$, see
(\ref{rk}), which is not decreasing, and thus $\lim_{n\ra\infty}
V_n=1-1/L$ for any $W$. \carre

\acknowledgement

The work of Luc Pronzato and Henry P.\ Wynn has been supported in
part by the IST Programme of the European Community, under the
PASCAL Network of Excellence, IST-2002-506778. This publication
only reflects the authors' views.


\bibliographystyle{plain}
\bibliography{xampl,test}

\end{document}